\magnification 1200
\overfullrule=0pt

\font\tenmsb=msbm10
\font\sevenmsb=msbm7
\font\fivemsb=msbm5
\font\tenmsa=msam10
\font\sevenmsa=msam7
\font\fivemsa=msam5
\newfam\msafam
\newfam\msbfam
\def\bbb{\fam\msbfam}
\def\Bbb#1{{\bbb{#1}}}
\textfont\msbfam=\tenmsb 
\scriptfont\msbfam=\sevenmsb
\scriptscriptfont\msbfam=\fivemsb
\textfont\msafam=\tenmsa 
\scriptfont\msafam=\sevenmsa
\scriptscriptfont\msafam=\fivemsa
\def\hexnumber@#1{\ifcase#1 0\or1\or2\or3\or4\or5\or6\or7\or8\or9\or
	A\or B\or C\or D\or E\or F\fi }
\edef\msa@{\hexnumber@\msafam}
\edef\msb@{\hexnumber@\msbfam}

\font\teneuf=eufm10
\font\seveneuf=eufm7
\font\fiveeuf=eufm5
\newfam\euffam
\textfont\euffam=\teneuf 
\scriptfont\euffam=\seveneuf
\scriptscriptfont\euffam=\fiveeuf
\def\got{\fam\euffam}

\font\eightrm=cmr8
\font\eightbf=cmbx8
\font\eightit=cmti8
\font\eightsl=cmsl8
\font\eightmus=cmmi8
\def\small{\let\rm=\eightrm \let\bf=\eightbf \let\it=\eightit
\let\sl=\eightsl \let\mus=\eightmus \baselineskip=9.5pt minus .75pt
}

\catcode`\@=11
\font\tenmsa=msam10
\font\sevenmsa=msam7
\font\fivemsa=msam5
\font\tenmsb=msbm10
\font\sevenmsb=msbm7
\font\fivemsb=msbm5
\newfam\msafam
\newfam\msbfam
\textfont\msafam=\tenmsa  \scriptfont\msafam=\sevenmsa
\scriptscriptfont\msafam=\fivemsa
\textfont\msbfam=\tenmsb  \scriptfont\msbfam=\sevenmsb
\scriptscriptfont\msbfam=\fivemsb
\def\hexnumber@#1{\ifcase#1 0\or1\or2\or3\or4\or5\or6\or7\or8\or9\or
A\or B\or C\or D\or E\or F\fi }
\edef\msa@{\hexnumber@\msafam}
\edef\msb@{\hexnumber@\msbfam}
\mathchardef\ltimes="2\msb@6E 
\mathchardef\boxtimes="2\msa@02 
\mathchardef\twoheadrightarrow="3\msa@10 
\mathchardef\twoheadleftarrow="3\msa@11 
\def\Box{\hfill$\sqcap\hskip-6.5pt\sqcup$}

\def\eb{{\bf e}}

\def\tb{{\bf t}}
\def\ub{{\bf u}}
\def\vb{{\bf v}}
\def\wb{{\bf w}}

\def\Ab{{\bf A}}
\def\Eb{{\bf E}}
\def\Fb{{\bf F}}
\def\Hb{{\bf H}}
\def\Ib{{\bf I}}
\def\Tb{{\bf T}}
\def\Ub{{\bf U}}
\def\Rb{{\bf R}}
\def\Kb{{\bf K}}
\def\Mb{{\bf M}}
\def\Vb{{\bf V}}
\def\Wb{{\bf W}}

\def\AA{{\Bbb A}}

\def\CC{{\Bbb C}}
\def\KK{{\Bbb K}}
\def\GG{{\Bbb G}}

\def\NN{{\Bbb N}}
\def\PP{{\Bbb P}}

\def\ZZ{{\Bbb Z}}

\def\Ec{{\cal E}}
\def\Fc{{\cal F}}

\def\Hc{{\cal H}}
\def\Kc{{\cal K}}
\def\Lc{{\cal L}}

\def\Oc{{\cal O}}

\def\Gr{{\GG{\rm r}}}

\centerline {{\bf AFFINE QUANTUM GROUPS}}
\vskip 1mm
\centerline {{\bf AND EQUIVARIANT K-THEORY}}

\vskip .7cm

\centerline {\bf E. VASSEROT
\footnote*{\fiverm The author is partially supported
by EEC grant no. ERB FMRX-CT97-0100.\hfill}}

\vskip 1cm

\noindent{\bf Introduction.} The purpose of this paper is to give complete
proofs of the $K$-theoretic construction of the quantized 
enveloping algebra of the affine Lie algebra $\widehat{\got gl}(n)$
sketched in [GV]. This construction gives a geometric 
classification of finite dimensional simple modules of $\widehat{sl}(n)$
in terms of intersection cohomology of graded
nilpotent orbits of type $A$.
The principal results were obtained jointly with V. Ginzburg. 
The formula for the Drinfeld polynomials of the simple modules
in the last section was not given in [GV].
The basic construction is a rephrasing
of the $K$-theoretic construction of the Iwahori-Hecke algebra 
given by Ginzburg.
Since a detailed proof is available for the case of the 
Iwahori-Hecke
algebra (see [CG]), the present paper deals mainly with the parts of the proof
which are different from the affine Hecke algebra one, more precisely 
Theorem 1 whose proof is inspired by [V]
and Theorem 2 whose proof is inspired by [BLM].

\vskip2mm

\noindent The first section of the paper contains some background 
and some technical lemmas on the K-groups of various flag varieties. 
We give complete proofs since we were not able to find a proper reference.
Similar results in homology may be found in [CG, Chapter 2].
The K-theoretic construction of the quantized enveloping algebra is done in 
Section two (Theorem 1 and Theorem 2). 
The simple and standard modules are studied in the third section. 
The reader should be warned that we give a geometric construction of the quantum
group of type $\widehat{\got gl}(n)$ while [GV] deals with type
$\widehat{\got sl}(n)$. 

\vskip3mm

\noindent {\bf 1. The convolution algebra in equivariant K-theory.}

\vskip3mm

\noindent {\bf 1.0. Algebraic equivariant K-theory.}
Throughout the paper a {\sl variety} means a complex quasi-projective
algebraic variety and {\sl sub-variety} means closed
subset for the Zariski topology. Given a complex linear algebraic 
group $G$ and a $G$-variety $X$, 
let ${\rm Vect}^G(X)$ (resp. ${\rm Coh}^G(X)$) denote
the set of isomorphism classes of 
$G$-equivariant vector bundles (resp. $G$-equivariant coherent sheaves) on $X$.
Let $K^G(X)$ be the complexified Grothendieck group of ${\rm Coh}^G(X)$.
The $K$-group has a natural $R(G)$-module structure 
where $R(G)$ is the complexified representation ring of $G$.
Given $\Fc\in{\rm Coh}^G(X)$ let $[\Fc]$ denote its class in $K^G(X)$. 
Recall a few properties of the equivariant K-theory (see [CG] for
more details).

\vskip2mm

\noindent (a) For any proper map $f\,:\,X\to Y$ between two $G$-varieties 
$X$ and $Y$ there is a derived direct image $R f_*\,:\,K^G(X)\to K^G(Y)$. 
The map $R f_*$ is a group homomorphism.

\vskip2mm

\noindent (b) If $f\,:\,X\to Y$ is flat (for instance an open embedding)
or is a closed embedding of a smooth $G$-variety, and $Y$ is smooth,
there is an inverse image $f^*\,:\,K^G(Y)\to K^G(X)$. The map $f^*$ is
a ring homomorphism for the product defined in subsection (g).

\vskip2mm

\noindent (c) Given $G$-varieties $X_1,X_2,Y_1,Y_2$ and a Cartesian square
$$\matrix{
X_1&\buildrel f_1\over\longrightarrow&Y_1\cr
{\scriptstyle g}\downarrow&&\downarrow{\scriptstyle h}\cr
X_2&\buildrel f_2\over\longrightarrow&Y_2\cr
}$$
such that $f_1, f_2$ are proper and $g,h$ are flat we have 
$Rf_{1\,*}\,g^*\,=\,h^*\,Rf_{2\,*}$ (proper base change).

\vskip2mm

\noindent (d) Fix a smooth $G$-variety $X$ and consider
a smooth $G$-subvariety $i\,:\,Y\hookrightarrow X$. Let 
$T_Y^*X\in{\rm Vect}^G(Y)$
be the conormal bundle of $Y$ in $X$. Set
$$\Lambda(T_Y^*X)=\sum_i(-1)^i[\Lambda^iT_Y^*X]\in K^G(Y).$$
Then $i^*\,Ri_*\,:\,K^G(Y)\to K^G(Y)$ is the product
by $\Lambda(T_Y^*X)$. 

\vskip3mm

\noindent (e) Given a $G$-equivariant vector bundle 
$\pi_{_E}\,:\,E\to X$ 
with zero section $\sigma_{_E}\,:\,X\hookrightarrow E$ 
over a smooth $G$-variety $X$,
the pull-back map $\pi^*_{_E}\,:\,K^G(X)\to K^G(E)$ is invertible
(Thom Isomorphism Theorem) and $\sigma_{_E}^*={\pi_{_E}^*}^{^{-1}}$.

\vskip2mm

\noindent (f) Consider a smooth $G$-variety $X$. 
For any $g\in G$ the $g$-fixed point subvariety, denoted $X^g$, is smooth. 
Let $i\,:\,X^g\hookrightarrow X$ be the inclusion. 
Suppose moreover that $G$ is Abelian. Then $X^g$ is a $G$-variety.
The complexified Grothendieck ring $R(G)$ is identified with the ring of
regular functions on $G$ and $K^G(X)$, $K^G(X^g)$,
can be viewed as sheaves over $G$.
Then, the direct image $Ri_*$ induces an isomorphism of the 
localized K-groups $K^G(X^g)_g{\buildrel\sim\over\longrightarrow}K^G(X)_g$
and $\Lambda(T^*_{X^g}X)$ is invertible in $K^G(X^g)_g$.
Moreover if $X$ is complete and $p\,:\,X\to pt$ and
$q\,:\,X^g\to pt$ are the (proper) projections, then
for any sheaf $\Fc\in{\rm Coh}^G(X)$ the Lefschetz formula holds :
$$Rp_*[\Fc]=Rq_*\biggl(\Lambda(T^*_{X^g}X)^{-1}\otimes [\Fc]\biggr).$$

\vskip3mm

\noindent (g) If $X$ is a smooth $G$-variety, there is a derived tensor 
product ${\buildrel L\over\otimes}$ making $K^G(X)$ a ring. Given 
$\Fc_1,\Fc_2\in{\rm Coh}^G(X)$ the product is 
$[\Fc_1]\,{\buildrel L\over\otimes}\,[\Fc_2]=
\delta^*[\Fc_1\boxtimes\Fc_2],$
where $\boxtimes$ is the external tensor product and
$\delta\,:\,X\hookrightarrow X\times X$ is the diagonal embedding.

\vskip2mm

\noindent (h) Given two $G$-subvarieties $Z_1,Z_2\subseteq X$, the derived 
tensor product on $X$ induces a map 
$K^G(Z_1)\otimes K^G(Z_2)\to K^G(Z_1\cap Z_2)$.
This map depends of the ambient smooth variety $X$.
In particular,

\vskip3mm

\noindent{\bf Lemma 1.} {\sl Given two smooth
$G$-subvarieties $Z_1,Z_2\subset X$
with conormal bundles $T^*_{Z_1}X$, $T^*_{Z_2}X$, set 
$Z=Z_1\cap Z_2$ and 
$N={T^*_{Z_1}X|}_{Z}\cap{T^*_{Z_2}X|}_{Z}\in {\rm Vect}^G(Z)$.
If $Z$ is smooth and
${TZ_1|}_{Z}\cap{TZ_2|}_{Z}=TZ$, then for any 
$\Fc_1\in{\rm Vect}^G(Z_1), \Fc_2\in{\rm Vect}^G(Z_2)$ 
$$[\Fc_1]\,{\buildrel L\over\otimes}\,[\Fc_2]=
\sum_i(-1)^i[\Lambda^i N\otimes{\Fc_1|}_{Z}\otimes{\Fc_2|}_{Z}]
\in K^G(Z).$$
}

\noindent{\sl Proof.} 
Consider the deformation to the normal bundle $T_ZX$ of $Z$ in $X$ :
we get a flat family $M_ZX\to\CC$ whose fiber over $z\in\CC$ is $T_ZX$ if
$z=0$ and $X$ else. In particular $M_ZZ_1$ and $M_ZZ_2$ can be viewed 
as closed subschemes in $M_ZX$. Thus, the intersection product can be computed
on the special fiber (over $z=0$). We are reduced to the following situation : 
$X$ is a vector space with a linear $G$-action and $Z_1,Z_2$ are subspaces
stable by the action of $G$.
The result follows by using Koszul resolutions.
\Box

\vskip3mm

\noindent {\bf 1.1. Composition of correspondences and convolution product.}
Let $G$ be a linear complex algebraic group and $M_1$, $M_2$ and $M_3$
be smooth $G$-varieties. Let 
$$q_{_{ij}}\,:\,M_1\times M_2\times M_3
\to M_i\times M_j$$
be the projection along the factor not named.
The $G$-action on each factor induces a natural $G$-action on the Cartesian
product such that the projections $q_{_{ij}}$ are $G$-equivariant. Let
$Z_{_{12}}\i M_1\times M_2$ and $Z_{_{23}}\i M_2\times M_3$
be $G$-stable closed subvarieties. Assume that the restriction of $q_{_{13}}$
to $q_{_{12}}^{^{-1}}Z_{_{12}}\cap q_{_{23}}^{^{-1}}Z_{_{23}}$ is proper
and let $Z_{_{12}}\circ Z_{_{23}}$ be its image. The $G$-subvariety
$Z_{_{12}}\circ Z_{_{23}}$ of $M_1\times M_3$
is called the {\sl composition} of $Z_{_{12}}$ and $Z_{_{23}}$.
Define a {\sl convolution} map 
$$\star\,:\,K^G(Z_{_{12}})\otimes K^G(Z_{_{23}})\to 
K^G(Z_{_{12}}\circ Z_{_{23}})$$
as follows. Let $\Fc_{_{12}}$, $\Fc_{_{23}}$
be two equivariant coherent sheaves respectivelly on
$Z_{_{12}}$ and $Z_{_{23}}$. Set
$$[\Fc_{_{12}}]\star[\Fc_{_{23}}]=
R q_{_{13*}}\biggl(q^*_{_{12}}[\Fc_{_{12}}]
\buildrel {L}\over\otimes q^*_{_{23}}[\Fc_{_{23}}]\biggr).$$
In this formula, the upper star stands for the pullback morphism, well-defined
on smooth maps. 

\vskip3mm

\noindent Let us now recall some technical results
in intersection theory which will be used 
in the proof of surjectivity theorem (Theorem 2). 
Suppose that $F_i$ is a smooth $G$-variety for $i=1,2,3$. 
Let $G$ act on $M_i=T^*F_i$ in such a way that
the projection $\pi_i\,:\,M_i\longrightarrow F_i$ is $G$-equivariant. 
Let 
$$p_{_{ij}}\,:\,F_1\times F_2\times F_3\to F_i\times F_j$$
be the obvious projection.
Fix a locally closed smooth $G$-subset
$O_{_{ij}}\i F_i\times F_j$ for all $(ij)=(12),(23),(13)$.
Denote by $Z_{_{ij}}=T^*_{O_{ij}}(F_i\times F_j)$ the conormal bundle
and let $\pi_{ij}\,:\,Z_{_{ij}}\to O_{_{ij}}$ be the projection.
Set $Z=q^{^{-1}}_{_{12}}Z_{_{12}}\cap q^{^{-1}}_{_{23}}Z_{_{23}}$
and $O=p^{^{-1}}_{_{12}}O_{_{12}}\cap p^{^{-1}}_{_{23}}O_{_{23}}$ .
Suppose that 

\vskip1mm

\noindent{(a)} the map $q_{_{13}}|_Z$ is proper,

\vskip1mm

\noindent{(b)} the projections 
$O_{_{12}}\to F_2$ and $O_{_{23}}\to F_2$ 
are smooth fibrations,
\vskip1mm

\noindent{(c)} 
the restriction $p=p_{_{13}}|_O$ is a smooth and 
proper fibration over $O_{_{13}}$. 

\vskip1mm

\noindent As usual, $Tp$ stands for the relative tangent sheaf along the fibers
of $p$.

\vskip3mm

\noindent{\bf Proposition 2.} {\sl 
If the hypothesis (a-c) hold 
then, $Z_{_{12}}\circ Z_{_{23}}=Z_{_{13}}$.
If moreover the $O_{_{ij}}$ are closed, then 
$\forall\Fc_{_{12}}\in{\rm Vect}^G(O_{_{12}})$,
$\forall\Fc_{_{23}}\in{\rm Vect}^G(O_{_{23}})$, 
$$\pi_{_{12}}^*[\Fc_{_{12}}]\star\pi_{_{23}}^*[\Fc_{_{23}}]=\sum_i(-1)^i
\pi^*_{_{13}}R p_*
[\Lambda^iT p\otimes 
p^*_{_{12}}\Fc_{_{12}}|_O\otimes p^*_{_{23}}\Fc_{_{23}}|_O].$$
}

\vskip1mm

\noindent{\sl Proof.} 
The property (a) implies that 
$$Z\cap (F_1\times T^*F_2\times F_3)=
(T^*_{O_{12}}(F_1\times F_2)\times F_3)\cap 
(F_1\times T^*_{O_{23}}(F_2\times F_3))=O,$$
where $O$ is identified with its image via the zero section.
Thus,
$$(TO_{_{12}}\times TF_3)+(TF_1\times TO_{_{23}})
=T(F_1\times F_2\times F_3).\leqno{\rm (d)}$$
On the other hand,
the hypothesis (b) implies that the projection $O\to F_2$
is a smooth fibration.
Thus, $p_{_{12}}^{^{-1}}O_{_{12}}$
and $p_{_{23}}^{^{-1}}O_{_{23}}$ 
intersect transversally along $O$, i.e. 
$$TO=(TO_{_{12}}\times TF_3)\cap (TF_1\times TO_{_{23}}).\leqno{\rm (e)}$$
The conormal bundle to $q_{_{12}}^{^{-1}}Z_{_{12}}$
in $M_1\times M_2\times M_3$, say $N_{_{12}}$, is 
the pull back of the vector bundle 
$$TO_{_{12}}\oplus T^*_{O_{12}}(F_1\times F_2)$$
by the obvious projection
$q_{_{12}}^{^{-1}}Z_{_{12}}\to O_{_{12}}.$
Let $\pi\,:\,Z\to O$ be the projection.
The formulas (d) and (e) give
$$Tp=(TO_{_{12}}\times F_3)\cap (F_1\times TO_{_{23}})
\quad{\rm and}\quad
{N_{_{12}}}|_Z\cap{N_{_{23}}}|_Z=\pi^*Tp
.$$
Since the projections $O_{_{12}}\to F_2$ and
$O_{_{23}}\to F_2$ are smooth fibrations,
we get smooth fibrations
$$T^*_{O_{12}}(F_1\times F_2)\subset(T^*F_1)\times F_2\to F_2
\quad{\rm and}\quad
T^*_{O_{23}}(F_2\times F_3)\subset F_2\times(T^*F_3)\to F_2,$$
and $Z$ is isomorphic to
$$Z\simeq T^*_{O_{12}}(F_1\times F_2)\times_{F_2}
T^*_{O_{23}}(F_2\times F_3).$$
In particular, $Z$ is smooth. Moreover, since 
$p_{_{12}}^{-1}O_{_{12}}$ and $p_{_{23}}^{-1}O_{_{23}}$
intersect transversally (see (e)),
$q_{_{12}}^{-1}Z_{_{12}}$ and $q_{_{23}}^{-1}Z_{_{23}}$
intersect transversally also. Thus Lemma 1 implies that 
$$q^*_{_{12}}\pi_{_{12}}^*[\Fc_{_{12}}]
\buildrel L\over\otimes 
q^*_{_{23}}\pi_{_{23}}^*[\Fc_{_{23}}]
= \sum_i(-1)^i\pi^*
[\Lambda^iTp\otimes p^*_{_{12}}\Fc_{_{12}}|_O\otimes
p^*_{_{23}}\Fc_{_{23}}|_O].$$
Now $O\simeq O_{12}\times_{F_2}O_{23}$ and thus
$$T^*_O(F_1\times F_2\times F_3)\simeq T^*_{O_{12}}(F_1\times F_2)\times_{F_2}
T^*_{O_{23}}(F_2\times F_3)\simeq Z.$$
Since the map $p$ is a smooth fibration,
$$q_{_{13}}(T^*_O(F_1\times F_2\times F_3))=T^*_{O_{13}}(F_1\times F_3).$$
Thus $Z_{_{12}}\circ Z_{_{23}}=Z_{_{13}}$.
Finally consider the commutative square
$$\matrix{
O&\buildrel p\over\longrightarrow &O_{_{13}}\cr
\pi\uparrow&&\uparrow\pi_{_{13}}\cr
Z&\buildrel q_{_{13}}\over\longrightarrow &
Z_{_{13}}.
}$$
Since the vertical maps are vector bundles and
$q_{_{13}}|_Z$ is proper and surjective,
the square is Cartesian. The result follows by
smooth base change.
\Box

\vskip3mm

\noindent As usual, for any locally-closed subset $X$
let $\overline X$ be the Zariski closure and set
$\partial X=\overline X\setminus X$.
Set $O_T=O\cap p_{_{13}}^{^{-1}}O_{_{13}}$.

\vskip3mm

\noindent{\bf Corollary 3.} {\sl Suppose in addition 
to {\rm (a-c)} that

\vskip1mm

\item{\rm (f)} $O_{_{13}}\subset O_{_{12}}\circ O_{_{23}}$ is open, 

\vskip1mm

\item{\rm (g)} $O_{_{13}}\cap(\partial O_{_{12}}\circ\overline O_{_{23}})=
O_{_{13}}\cap(\overline O_{_{12}}\circ\partial O_{_{23}})=\emptyset$,

\vskip1mm

\noindent 
Then $Z_{_{13}}\subset\overline Z_{_{12}}\circ \overline Z_{_{23}}$
is open and, for all $\Fc_{_{12}}\in{\rm Vect}^G(\overline Z_{_{12}})$,
$\Fc_{_{23}}\in{\rm Vect}^G(\overline Z_{_{23}})$, 
the restriction of $[\Fc_{_{12}}]\star[\Fc_{_{23}}]$ 
to $Z_{_{13}}$ is 
$$\sum_i(-1)^i\pi^*_{_{13}}Rp_*[
\Lambda^iTp\otimes
p^*_{_{12}}(\Fc_{_{12}}|_{O_{_{12}}})\otimes
p^*_{_{23}}(\Fc_{_{23}}|_{O_{_{23}}})].$$
}

\vskip1mm

\noindent{\sl Proof.} 
Set $Z_T=Z\cap q_{_{13}}^{^{-1}}Z_{_{13}}$.
The hypothesis {\rm (g)} implies that the commutative square
$$\matrix{
q_{_{12}}^{^{-1}}\overline Z_{_{12}}\cap q_{_{23}}^{^{-1}}\overline 
Z_{_{23}}&
\buildrel q_{_{13}}\over\longrightarrow 
&\overline Z_{_{12}}\circ\overline Z_{_{23}}\cr
\uparrow&&\uparrow\cr
Z_T&\buildrel q_{_{13}}\over\longrightarrow & Z_{_{13}} }$$
is Cartesian. Moreover the vertical maps are open because of {\rm (f)}.
Thus the restriction of $[\Fc_{_{12}}]\star[\Fc_{_{23}}]$ 
to $Z_{_{13}}$ is 
$$Rq_{_{13\,*}}\bigr( 
(q^*_{_{12}}[\Fc_{_{12}}]\buildrel L\over\otimes q^*_{_{23}}[\Fc_{_{23}}])
|_{Z_T}\bigr).$$
Then, apply Proposition 2.
The square
$$\matrix{
O_T&\buildrel p\over\longrightarrow &O_{_{13}}\cr
\uparrow&&\uparrow\pi_{_{13}}\cr
Z_T&\buildrel q_{_{13}}\over\longrightarrow &
Z_{_{13}} }$$
is Cartesian (see the proof of Proposition 2)
and we are done.
\Box

\vskip3mm

\noindent{\bf Corollary 4.} {\sl  
Suppose that $F_3=pt$ and 
$O\subset F_1\times F_2$ is a smooth $G$-subvariety with 
conormal bundle $Z_O\subset M_1\times M_2$ such that

\vskip1mm

\item{\rm (h)} the projection $Z_O\to M_1$ is proper, 

\vskip1mm

\item{\rm (i)} both projections 
$p_{1,O}\,:\,O\to F_1$ and 
$p_{2,O}\,:\,O\to F_2$ 
are smooth fibrations and $p_{1,O}$ is proper.

\vskip1mm

\noindent If $\Fc\in{\rm Vect}^G(F_2)$ then,
$$[\Oc_{Z_O}]\star \pi_2^*[\Fc]=\sum_i(-1)^i\pi_1^*
Rp_{1,O\,*}
[\Lambda^iTp_{1,O}\otimes p_{2,O}^*\Fc]
\in K^G(M_1).$$
}

\vskip3mm

\noindent{\sl Proof.} 
It suffices to apply Proposition 2 with $O_{_{12}}=O$,
$\Fc_{_{12}}=\Oc_{O}$, $O_{_{23}}=F_2$
and $\Fc_{_{23}}=\Fc$ (see either [V, lemme 5]).
\Box

\vskip3mm

\noindent {\bf 1.2. Partitions, matrices, and 3-arrays.}
For any positive integer $n$ set $[n]=\{1,2,\cdots,n\}$.
Fix two integers $n,d\geq 1$.
Let ${\bf V}\subset\NN^n$ be the set of
all compositions ${\bf v} = (v_1,\ldots,v_n)$
of $d$ into $n$ non negative integers. Set 
${\bf e}_i=(\delta_{1i},\delta_{2i},...,\delta_{ni})$, so that
${\bf v}=\sum_iv_i{\bf e}_i.$
For ${\bf v}\in {\bf V}$ let $[{\bf v}]_i \i [d]$
be the  $i$-th segment of the composition {\bf v} and
put $\overline v_i=v_1+v_2+...+v_i$.
Thus $[\vb]_i=[1+\overline v_{i-1},\overline v_i]$
and $[\vb]=([\vb]_1,[\vb]_2,...,[\vb]_n)$ is a partition
of $[d]$ into $n$ subsets. 

\vskip3mm

\noindent Let {\bf M} be the set of $(n\times n)$-matrices
$A=(a_{ij})$ with non-negative integral
entries such that $\sum_{i,j} a_{ij} = d$.
For  $\vb, \wb \in {\bf V}$ we put 
$${\bf M}(\vb,\wb) = \biggl\{ A\in {\bf M}\,\biggl|\,
\sum_ja_{ij}=v_i , \,\sum_ia_{ij}=w_j\biggl\},$$
so that ${\bf M}  = \coprod_{\vb,\wb\in {\bf V}}
{\bf M}(\vb,\wb)$. 
For any $\vb,\wb$ we have a natural map
$${\bf m} = {\bf m}_{\vb,\wb}: S_d \longrightarrow 
{\bf M}(\vb,\wb)$$
which identifies ${\bf M}(\vb,\wb)$ with the double coset
space $S_\vb\backslash S_d/S_\wb$. Namely, ${\bf m}$
assigns to a permutation $\sigma\in S_d$ the matrix ${\bf m}^\sigma$
such that 
$${\bf m}^\sigma_{ij} = \sharp\bigl\{\alpha\in [\vb]_i\,|\,
\, \sigma(\alpha)\in [\wb]_j\bigl\}.$$
We introduce a partial order $\preceq$ on ${\bf M}(\vb,\wb)$
as follows. For $A = (a_{ij})$
and $B = (b_{ij})$ in ${\bf M}(\vb,\wb)$ we say that
$A \preceq B$, if for any $1\leq i<j\leq n$ we have
$$\sum_{_{r\leq i\,;\,s\geq j}}a_{rs}\leq\sum_{_{r\leq i\,;\,s\geq j}}b_{rs}$$
and for any $1\leq j<i\leq n$ we have
$$\sum_{_{r\geq i\,;\,s\leq j}}a_{rs}\leq\sum_{_{r\geq i\,;\,s\leq j}}b_{rs}.$$
Denote by $\leq$ the Bruhat order on $S_d$. The map
${\bf m}: (S_d,\leq)\rightarrow ({\bf M}(\vb,\wb),\preceq)$
is monotone (see [BLM , Lemma 3.6]), i.e.
$$\sigma\leq \tau\Rightarrow {\bf m}^\sigma\preceq {\bf m}^\tau,\qquad
\forall\sigma,\tau\in S_d.$$
The order $\preceq$ on ${\bf M} = \coprod {\bf M}(\vb,\wb)$
is the disjoint union of the orders on the components 
${\bf M}(\vb,\wb)$, i.e. elements of different
components are set to be incomparable. 
To a matrix $A$ we associate a partition $[A]$
of the set $[d]$ labelled by the set $[n]\times[n]$ as follows :
to a pair $(i,j)$ corresponds the segment 
$$[A]_{ij}=\biggl[1+\sum_{_{(h,k)<(i,j)}}a_{hk},\,
\sum_{_{(h,k)\leq(i,j)}}a_{hk}\biggr],$$
where $\leq$ is the right lexicographic order defined by 
$$(i,j)>(h,k)\,\Leftrightarrow\,(j>k\quad{\rm or}\quad(j=k\quad{\rm
and}\quad i>h)).$$
similarly let $A_1,A_2$ be partitions of $[d]$ such that
$$A_{1i}=\bigcup_j[A]_{ij}\qquad{\rm and}\qquad
A_{2j}=\bigcup_i[A]_{ij}.$$

\vskip3mm

\noindent Let {\bf T} be the set of 3-arrays
$(t_{ijk})_{1\leq i,j,k\leq n}$ of non-negative integers,  such that 
$\sum_{i,j,k} t_{ijk} = d$. 
For any $\ub,\vb,\wb\in {\bf V}$ put
$${\bf T}(\ub,\vb,\wb) = 
\biggl\{  T\in {\bf T}\,\biggl|\,\sum_{j,k} t_{ijk} = u_i,\,\,
\sum_{i,k} t_{ijk} = v_j, \,\, \sum_{i,j} t_{ijk} = w_k\biggl\}.$$
If $1\leq i\leq j \leq 3$ and $T\in {\bf T}$,
let $T_{ij}\in\Mb$ be the matrix
obtained by summing the entries of $T$ with respect to the indices not named.
Given two matrices $A\in{\bf M}(\tb,\ub)$ and 
$B\in{\bf M}(\vb,\wb)$ put 
$${\bf T}(A,B)=\{T\in{\bf T}\,|\,T_{12}=A, T_{23}= B\},\qquad
{\bf M}(A,B)=\{T_{13}\,| \,T\in\Tb(A,B)\}.$$
The set ${\bf T}(A,B)$ is empty unless $\ub=\vb$. 

\vskip3mm

\noindent For any partition $I=(I_1,I_2,...,I_n)$ of the set $[d]$ into $n$ 
subsets let
$$S_I=S_{I_1}\times S_{I_2}\times\cdots\times S_{I_n}$$
be the subgroup of $S_d$ consisting of permutations which preserve
each subset. In particular, if $\vb\in\Vb$ and $A\in\Mb$ put
$S_\vb = S_{[\vb]}$ and $S_A = S_{[A]}.$

\vskip3mm

\noindent {\bf 1.3. The Flag manifold.} Fix $d,n \geq 1$. 
The complex linear group of rank $d$ is denoted by $GL_d$. 
Let $(\epsilon_1,\epsilon_2,...,\epsilon_d)$ be
the canonical basis of $\CC^d$.
Let $F$ be the variety of $n$-step partial flags in $\CC^d$, 
i.e. filtrations of vector spaces
$$D = \{0=D_0\subseteq D_1\subseteq D_2 \subseteq ...
\subseteq D_n =\CC^d\}.$$ 
The connected components of $F$
are parametrized by compositions ${\bf v} = (v_1, ..., v_n)$ of $d$ : 
the component $F_{\bf v}$ consists of
flags $D$ such that $\dim D_i/D_{i-1} = v_i$. For any composition let
$D_{\vb}\in F_{\vb}$ be the flag such that
$$D_{\vb,i}=\bigoplus_{j\in[1,\overline v_i]}\CC\epsilon_j.$$
Denote by $P_\vb$ the isotropy subgroup of $D_\vb$ in $GL_d$.
The variety $F_\vb$ is thus identified with the quotient $GL_d/P_\vb$.
Consider the diagonal action of $GL_d$ on $F\times F$ 
and $F\times F\times F$. The following result is essentially
stated in [BLM, Lemma 3.7].

\vskip3mm

\noindent{\bf Proposition 5. } {\it
(a) Given two partitions $\vb,\wb\in\Vb$,
the orbits of $GL_d$ in $F_{_\vb}\times F_{_\wb}$ are
parametrized by matrices $A\in {\bf M}(\vb,\wb)$ :
the orbit $O_A$ corresponding to
$A$ consists of pairs of flags $(D, D')$ such that
$\dim (D_i\cap D'_j) = \sum_{_{h\leq i,\,k\leq j}}a_{hk}$.
The closure of $O_A$ contains $O_B$ if and only if $B\leq A$
with respect to the Bruhat order on {\bf M}.  \hfill\break
(b) For any $A,B\in{\bf M}$ the set 
$p_{_{13}}(p_{_{12}}^{^{-1}}O_{_A} \cap p_{_{23}}^{^{-1}}O_{_B})\subset F^{^2}$
is stable with respect to the diagonal $GL_{_d}$ action and 
$$p_{_{13}}(p_{_{12}}^{^{-1}}O_{_A} \cap p_{_{23}}^{^{-1}}O_{_B})= 
\bigcup_{_{C\in\Mb(A,B)}} O_{_C}.$$
(c) Moreover, if $A,B\in{\bf M}$ there exists a unique
$A\circ B\in\Mb(A,B)$ such that 
$\Mb(A,B)\subset\{C\,|\,C\leq A\circ B\}.$\hfill\break
(d) More precisely, given $A,A',B,B'\in{\bf M}$ such that 
$A'\leq A$ and $B'\leq B$, 
$$(A',B')\neq (A,B)\,\Rightarrow \, 
\Mb(A',B')\subset\{C\,|\,C< A\circ B\}.$$
}

\vskip3mm

\noindent{\sl Proof.} Claims $(a)$ and $(b)$ are immediate.  
As for the part $(c)$ it suffices to prove that 
$p_{_{13}}(p_{_{12}}^{^{-1}}O_{_A} \cap p_{_{23}}^{^{-1}}O_{_B})$
is irreducible. But 
$p_{_{12}}\,:\,p_{_{12}}^{^{-1}}O_{_A} 
\cap p_{_{23}}^{^{-1}}O_{_B}\rightarrow O_{_A}$
is a fibration with an irreducible base and fiber. Thus 
$p_{_{12}}^{^{-1}}O_{_A} \cap p_{_{23}}^{^{-1}}O_{_B}$
is irreducible and so is its image under $p_{_{13}}$. 
Part $(d)$ follows easily.
\Box

\vskip3mm

\noindent For any $i=1,2,3$ let $p_{_i}\,:\,F^2\to F$ be the projections
along the component not named. For any $A$ 
let $p_{_{i,A}}$ be the restriction of $p_{_i}$ to $O_{_A}$.
This map is a smooth fibration but is not proper in general. 
Given a positive integer $a<d$ and a partition
${\bf v}$ of $d-a$ into $n$ summands and $i\neq j$ set
$$E_{_{ij}}({\bf v},a)={\rm diag}({\bf v})+a\cdot E_{ij}
\in{\bf M}({\bf v}+a\cdot{\bf e}_i,{\bf v}+a\cdot{\bf e}_j)$$
where 'diag' stands for the diagonal matrix with the prescribed eigenvalues
and $E_{ij}$ is the standard $n\times n$-matrix unit (1 at the spot $(i,j)$
and 0 elsewhere). If $j=i\pm 1$ the corresponding $GL_d$-orbits in $F^2$ are
closed and have the following description
$$O_{E_{_{i\,i+1}}({\bf v},a)}=\{(D,D')\in F^2\,|\,D'_k\subseteq D_k\quad
{\rm and}\quad {\rm dim}(D_k/D'_k)=a\,\delta_{ki}\},$$
$$O_{E_{_{i+1\,i}}({\bf v},a)}=\{(D,D')\in F^2\,|\,D_k\subseteq D'_k\quad
{\rm and}\quad {\rm dim}(D'_k/D_k)=a\,\delta_{ki}\}.$$
Thus both projections
$$F_{{\bf v}+a\cdot{\bf e}_i}
\twoheadleftarrow O_{E_{i\,i\pm 1}({\bf v},a)}
\twoheadrightarrow F_{{\bf v}+a\cdot{\bf e}_{i\pm 1}}$$ are proper 
with fibers isomorphic respectively 
to the Grassmannians $\Gr_a(a+v_i)$ and $\Gr_a(a+v_{i\pm 1})$. 

\vskip3mm

\noindent Given a 3-array $T\in\Tb$, put  
$$O_{_T}=p_{_{12}}^{^{-1}}O_{_{T_{12}}}\cap p_{_{23}}^{^{-1}}O_{_{T_{23}}}\cap
p_{_{13}}^{^{-1}}O_{_{T_{13}}}\subset F^3.$$
The set $O_{_T}$ is a $GL_d$-variety and may contain several 
$GL_d$-orbits in general.
For any $1\leq i<j\leq 3$ let $p_{_{ij}}\,:\,F^3\to F^2$ be the projection
along the component not named and let $p_{_{ij,T}}$ denote the restriction
of $p_{_{ij}}$ to $O_{_T}\subset F^3$.

\vskip3mm

\noindent Set
$\Rb=\CC[x_{_1}^{^{\pm 1}},x_{_2}^{^{\pm 1}},...,x_{_d}^{^{\pm 1}}],$ 
and let 
$\CC^{\times (d)}=(\CC^\times)^d/S_{_d}$ 
denote the $d$-fold 
symmetric product of $\CC^\times$.
For any partition $I$ of the set $[d]$ put 
$$\CC^{\times(I)}=(\CC^\times)^d/S_I.$$
The rings $\Rb^{(d)}=\Rb^{^{S_d}}$ and $\Rb^{(I)}=\Rb^{^{S_I}}$
may (and will) naturally be identified with the rings of regular functions on 
$\CC^{\times (d)}$ and $\CC^{\times (I)}$.
In order to simplify notation,
for $\vb\in\Vb$ and $A\in{\bf M}$ set 
$$\matrix{
\CC^{\times (\vb)}=\CC^{\times ([\vb])},\quad\hfill&
\CC^{\times (A)}=\CC^{\times ([A])},\quad\hfill\cr\cr
\Rb^{(\vb)}=\Rb^{([\vb])},\quad\hfill&
\Rb^{(A)}=\Rb^{([A])}.\quad\hfill
}$$
If $J$ is another partition of $[d]$
consider the symmetrizer ${\got S}_I^J$ and
the projection $\rho_I^J$ such that
$$\matrix{
{\got S}_I^J\,:\,
\Rb^{^{S_I\cap S_J}}\to\Rb^{^{(J)}},\qquad
f\mapsto\sum_{\sigma\in S_J/S_I\cap S_J}
\sigma(f),\hfill\cr\cr
\rho_I^J\,:\,
(\CC^\times)^d/S_I\cap S_J\to\CC^{\times (J)}.\hfill
}$$
The pull-back map 
$\rho_I^{J\,*}\,:\,\Rb^{^{(J)}}\to\Rb^{^{S_I\cap S_J}}$
is the inclusion of rings.
As before, if $\vb,\wb\in\Vb$ we write ${\got S}_\vb^\wb$ and $\rho_\vb^\wb$
instead of ${\got S}_{_{[\vb]}}^{^{[\wb]}}$
and $\rho_{_{[\vb]}}^{^{[\wb]}}$. 

\vskip3mm

\noindent {\bf Proposition 6.} {\sl 
Fix $\vb,\vb_1,\vb_2\in\Vb$, $A\in{\bf M}(\vb_1,\vb_2)$ and $i=1,2$. 

\noindent(a) 
$K^{GL_d}(F_{_\vb})\simeq\Rb^{^{(\vb)}}$  and
$K^{GL_d}(O_{_A})\simeq\Rb^{^{(A)}}$ as $\CC$-algebras. 

\noindent(b) The map $p_{_{i,A}}\,:\,O_{_A}\to F_{_{\vb_i}}$ 
is a smooth fibration. The inverse image morphism $p_{_{i,A}}^*$
in equivariant $K$-theory is identified with the pull-back map 
$\rho_A^{A_i\,*}\,:\,\Rb^{^{(A_i)}}\to\Rb^{^{(A)}}$.
If $O_{_A}$ is closed then the direct image morphism $Rp_{_{i,A\,*}}$ is 
well defined and 
$$Rp_{_{i,A\,*}}(f)={\got S}_A^{A_i}
\biggl(f\,\prod_{(s,t)}(1-x_s/x_t)^{-1}\biggr),$$
where the product ranges over all couples 
$$(s,t)\in [A]_{jl}\times [A]_{jk},\qquad{\rm if}\quad i=1,$$
$$(s,t)\in [A]_{lj}\times [A]_{kj},\qquad{\rm if}\quad i=2,$$
where $1\leq k<l\leq n,\quad j\in[n].$
}

\vskip3mm

\noindent {\sl Proof.} Fix $A\in\Mb(\vb,\wb)$.
Let $L_{_i}\subset\CC^{^d}$ (resp. $C_{_j}\subset\CC^{^d}$) be the 
subspace generated by the $\epsilon_t$ such that $t\in\bigcup_k [A]_{_{ik}}$ 
(resp. $t\in\bigcup_k [A]_{_{kj}}$).
Consider the flags $D\in F_{\vb}$ and $D'\in F_{\wb}$ such that 
$$D_{_i}=L_{_1}\oplus L_{_2}\oplus\cdots\oplus L_{_i}
\quad{\rm and}\quad
D'_{_j}=C_{_1}\oplus C_{_2}\oplus\cdots\oplus C_{_j}.$$
By definition $O_{_A}$ is the orbit of the pair 
$(D,D')$ for the diagonal action of $GL_{_d}$ on $F_{\vb}\times F_{\wb}$.
Let $P_D$ and $P_{D'}$ be the isotropy groups of $D$ and $D'$.
The induction property in equivariant $K$-theory gives an isomorphism
of the ring $K^{GL_d}(O_{_A})$ with the complexified 
representations ring of the group $P_D\cap P_{D'}$.
The groups $\prod_iGL(L_{_i})$ and $\prod_jGL(C_{_j})$ are 
isomorphic to Levi subgroups of $P_D$ and $P_{D'}$. 
Moreover $\prod_{i,j}GL(L_{_i}\cap C_{_j})$ is isomorphic to the
reductive part of $P_D\cap P_{D'}$. Part $(a)$ follows. By 
definition of induction, the map
$$R(P_D)
{\buildrel\sim\over\longrightarrow}
K^{GL_d}(F_{\vb})
{\buildrel p_{1,A}^*\over\longrightarrow} 
K^{GL_d}(O_A) 
{\buildrel\sim\over\longrightarrow}
R(P_D\cap P_{D'})$$
is the restriction map, i.e.
the pull-back by $\rho_A^{A_1}$. 
As for the direct image fix, for instance, $i=1$ and suppose that $O_A$ 
is closed. The fibers of the projection
$p_{_{1,A}}\,:\, O_A\to F_\vb$ are isomorphic to 
the product of partial flag manifolds
$$\prod_{i=1}^nF_{(a_{i1},a_{i2},...,a_{in})}.$$ 
Thus the class $T_A^*\in K^{GL_d}(O_A)$ of the relative cotangent bundle
to the fibers of $p_{1,A}$ is
$$T_A^*=\sum_{k<l\atop j}\sum_{s\in[A]_{jl}}\sum_{t\in[A]_{jk}}
x_s^{-1}\cdot x_t,$$
where the $k$-th elementary symmetric polynomial
in the $x_s$, $s\in[A]_{ij}$, stands for the class
in $K^{GL_d}(O_A)$ of the vector bundle whose
fiber at $(D,D')\in O_A$ is
$$\bigwedge^k\biggl(
{D_i\cap D'_j\over (D_{i-1}\cap D'_j)+(D_i\cap D'_{j-1})}
\biggr).$$
Set $\Lambda T_A^*=\sum_i(-1)^i\Lambda^iT_A^*$,
where $\Lambda^iT_A^*$ is the $i$-th wedge of $T_A^*$.
Then, the direct image by $p_{1,A}$ is the map 
$$K^{GL_d}(O_A)\longrightarrow K^{GL_d}(F_{A_1}),\qquad
[\Fc]\mapsto{\got S}_A^{A_1}\biggl([\Fc]\otimes\Lambda(T_A^*)^{-1}\biggr)$$
(see [V, Lemme 4]). Part $(b)$ follows.
\Box

\vskip3mm

\noindent{\bf 1.4. The Steinberg variety and the convolution algebra.}
Fix a formal variable $q$ and set $\AA=\CC[q,q^{-1}]$, $\KK=\CC(q)$.
For any complex vector space $V$ set $V_\AA=V\otimes_\CC\AA$ and
$V_\KK=V\otimes_\CC\KK$. Similarly, if $V$ is a $\AA$-module put
$V_\KK=V\otimes_\AA\KK$.
Set $G=GL_d\times\CC^\times$ and use the notations of the previous sections. 
In particular $F=\coprod_{\vb\in\Vb}F_\vb$ is the $n$-step flag variety. 
For any $A\in\Mb$
let $Z_A= T^*_{O_A} (F\times F)\i T^*F^2$ be the conormal bundle to $O_A$.
Set 
$$Z_{\vb,\wb}=\bigcup_{A\in\Mb(\vb,\wb)}\overline Z_A
\quad\i\quad T^*(F_{\vb}\times F_{\wb})
\quad{\rm and}\quad
Z=\bigcup_{A\in\Mb}\overline Z_A\quad\i\quad T^*F^2.$$
Using the well-known isomorphism 
$T^*F=\{(D,x)\in F\times{\rm End}(\CC^d)\,|\,x(D_{i+1})\subseteq D_i
\}$
the set $Z$ is identified with the following closed subvariety of $T^*F^2$ 
$$Z=\{(D,D',x,x')\in T^*F^2\,|\, x=x'\}.$$
The variety $Z$ is called the Steinberg variety.
It is a reducible variety whose irreducible components are the $\overline Z_A$,
$A\in\Mb$.
Let $q_{_{12}}$, $q_{_{23}}$ and $q_{_{13}}$ be the projections
$T^*F^3\to T^*F^2$ along the component not named. 
It is known that $Z$ satisfies the following two properties 

\vskip2mm

\item{\bf .} the restriction of the projection $q_{_{13}}$ to
$q_{_{12}}^{^{-1}}Z\cap q_{_{23}}^{{-1}}Z$  is proper,

\vskip2mm

\item{\bf .} the composition 
$Z\circ Z=q_{_{13}}(q_{_{12}}^{^{-1}}Z\cap q_{_{23}}^{{-1}}Z)$ 
is equal to $Z$.

\vskip3mm

\noindent 
The first property is obvious, the second one follows from Proposition 2
(see also [CG, Chapter 6]). The group $G$ acts naturally on $Z$ : the linear 
group
$GL_d$ acts diagonaly and $z\in\CC^\times$ acts by scalar multiplication
by $z^{-2}$ along the fibers. The convolution product 
$$\star\,:\,K^G(Z)\otimes K^G(Z)\to K^G(Z)$$
endows $K^G(Z)$ with the structure of an associative $\AA$-algebra.
For any $A\in\Mb$ set $Z_{_{\leq A}}=\bigcup_{B\leq A}Z_{_B}$.
The natural maps $K^G(Z_{\leq A})\to K^G(Z)$ induced by the closed embeddings
$Z_{\leq A}\hookrightarrow Z$ are injective and their images form a filtration
on the algebra $K^G(Z)$ indexed by ${\bf M}$. 
Proposition 5 implies that 
$Z_{_{\leq A}}\circ Z_{_{\leq B}}\subset Z_{_{\leq A\circ B}}$.
Thus, $K^G(Z_{_{\leq A}})\star K^G(Z_{_{\leq B}})\i K^G(Z_{_{\leq A\circ B}})$. 
On the other hand, the open immersion 
$Z_{_A}\hookrightarrow Z_{_{\leq A}}$ gives rise 
to a restriction map $K^G(Z_{\leq A})\to K^G(Z_A)$ 
and to a short exact sequence
$$0\to K^G(Z_{_{<A}})\to K^G(Z_{_{\leq A}})\to K^G(Z_{_A})\to 0,$$
where $Z_{_{<A}}$ stands for $\bigcup_{B<A}Z_{_{\leq B}}$
(see [CG, 5.5]). The chain of maps 
$$\,K^G(Z_{\leq A})\twoheadrightarrow
K^G(Z_{\leq A})/K^G(Z_{<A})
\buildrel\sim\over\rightarrow K^G(Z_A)
\buildrel\sim\over\rightarrow K^G(O_A)
\buildrel\sim\over\rightarrow \Rb^{^{(A)}}_{_\AA}$$
provides an identification of the associated graded of $K^G(Z)$ with
the $\AA$-module $\bigoplus_{_{A\in\Mb}}\Rb^{^{(A)}}_{\AA}$.
Since by Proposition 5 we have
$Z_{_{\leq A'}}\circ Z_{_{\leq B'}}\subset Z_{_{< A\circ B}}$
if $A'\leq A$, $B'\leq B$ and $(A',B')\neq (A,B)$, the convolution
product induces a product still denoted $\star$ on 
$\bigoplus_{_{A\in\Mb}}\Rb^{^{(A)}}_{\AA}$ such that
$$\Rb^{^{(A)}}_{_\AA}\star\Rb^{^{(B)}}_{_\AA}\subset\Rb^{^{(A\circ B)}}_{_\AA}.$$

\vskip1mm

\noindent{\bf Example.} Given two positive integers $a,b<d$,
and a partition $\vb$ of $d-a-b$ put $A=E_{_{i\,i+1}}(\vb+b\eb_i,a)$ 
and $B=E_{_{i\,i+1}}(\vb+a\eb_{i+1},b)$ (see Section 1.3). 
The orbits $O_{_A}$ and $O_{_B}$ are closed. 
The set $\Tb(A,B)$ reduces to a single element,
$T$, and $A\circ B=E_{_{i\, i+1}}(\vb,a+b)$. 
For simplicity take $n=2$, $i=1$. Then $d=a+b+v_1+v_2$,
$$A=\pmatrix{v_{_1}+b&a\cr 0&v_{_2}},\quad
B=\pmatrix{v_{_1}&b\cr 0&v_{_2}+a},\quad
A\circ B=\pmatrix{v_{_1}&a+b\cr 0&v_{_2}},$$
and the 3-array $T$ is defined by 
$$t_{121}=t_{212}=t_{211}=t_{221}=0\quad
{\rm and}\quad t_{122}=a,t_{112}=b,t_{111}=v_1,t_{222}=v_2.$$
Set $I_1=[1,v_1]$, $I_2=[1+v_1,b+v_1]$, $I_3=[1+b+v_1,d-v_2]$
and $I_4=[1+d-v_2,d]$. Thus,
$$S_A=S_{I_1\cup I_2}\times S_{I_3}\times S_{I_4},\quad
S_B=S_{I_1}\times S_{I_2}\times S_{I_3\cup I_4},\quad
S_{A\circ B}=S_{I_1}\times S_{I_2\cup I_3}\times S_{I_4}.$$
The restriction of $p_{_{13}}$ to 
$O_{_T}=p_{_{12}}^{^{-1}}O_{_A}\cap p_{_{23}}^{^{-1}}O_{_B}$
is a smooth and proper fibration whose fiber over 
$(D,D'')\in O_{_{A\circ B}}$ is
$$\{D'\in F\,|\,D_k\subseteq D'_k\subseteq D''_k\quad 
{\rm and\quad dim}(D'_k/D_k)=a\delta_{ki},\quad\forall k\}
\simeq\Gr_a(a+b).$$
Note that $q^{^{-1}}_{_{12}}Z_{_A}$ and $q^{^{-1}}_{_{23}}Z_{_B}$
are not in general position in $T^*F^3$.
Proposition 5 insures that we are in the
situation of Corollary 3.
Using the isomorphism above the convolution product 
$\star\,:\,K^G(Z_{_A})\otimes K^G(Z_{_B})\to K^G(Z_{_{A\circ B}})$ 
is identified with the map
$$\Rb^{^{(A)}}_{_\AA}\otimes\Rb^{^{(B)}}_{_\AA}\to\Rb^{^{(A\circ B)}}_{_\AA},
\qquad f\otimes g\mapsto {\got S}_{I_2\times I_3}^{I_2\cup I_3}
\biggl( fg\cdot\prod_{i\in I_2}
\prod_{j\in I_3}{1-q^2x_j/x_i\over 1-x_i/x_j}\biggr)$$
(see Proposition 6 and Corollary 3).
Now if $f\in\Rb^{(I_2\cup I_3)}$ then
$${\got S}_{I_2\times I_3}^{I_2\cup I_3}
\biggl(f\cdot\prod_{i\in I_2}\prod_{j\in I_3}{q^2x_j-x_i\over x_j-x_i}\biggr)=
f\cdot{\got S}_{I_2\times I_3}^{I_2\cup I_3}
\biggl(\prod_{i\in I_2}\prod_{j\in I_3}{q^2x_j-x_i\over x_j-x_i}\biggr)=
P(q)\cdot f$$
where $P(q)$ is a polynomial in $q$ which can be computed explicitely :
$$P(q)=q^{ab}{[a+b]!\over [a]!\,[b]!}\quad{\rm with}\quad
[k]!=[k]\,[k-1]\cdots[1],
\qquad[k]={q^{k}-q^{-k}\over q-q^{-1}}\qquad\forall k.$$
Observe that the zeros of $P$ are roots of unity.
Since
$$\prod_{i\in I_2}\prod_{j\in I_3}{q^2x_j-x_i\over x_j-x_i}=
\prod_{i\in I_2}\prod_{j\in I_3}(-x_i/x_j){1-q^2x_j/x_i\over 1-x_i/x_j},$$
and $\prod_{i\in I_2}x_i\in\Rb^{^{(B)}}_\AA$, 
$\prod_{j\in I_3}x_j\in\Rb^{^{(A)}}_\AA$ are invertible, 
the map
$\Rb^{^{(A)}}_{_\KK}\otimes\Rb^{^{(B)}}_{_\KK}\to\Rb^{^{(A\circ B)}}_{_\KK}$
is surjective, and remains surjective after specialization of $q$
to a non zero complex number distinct from a root of unity.

\vskip3mm

\noindent{\bf Proposition 7.} {\sl Fix a partition $\vb\in\Vb$ and
suppose that $A={\rm diag}\,(\vb)$ is a diagonal matrix.
Then $A\circ B=B$ and, given $f\in\Rb^{^{(A)}}_{_\AA}$,
$g\in\Rb^{^{(B)}}_{_\AA}$, we have
$$f\star g= f\, g,$$
where $f$ is identified with an element of $\Rb^{^{(B)}}_{_\AA}$
via the map $\rho_B^{B_1}$.
}

\vskip1mm

\noindent {\sl Proof.} The projection $p_{_{13}}$ restricts to an isomorphism
$p_{_{12}}^{^{-1}}O_{_A}\cap p_{_{23}}^{^{-1}}O_{_B}\to O_{_B}$.
Thus $A\circ B=B$ and $\Tb(A,B)=\{T\}$ where
$T$ is the array such that $t_{ijk} = \delta_{_{ij}} b_{jk}$.
The proposition follows from Propositions 5 and 6, and
Corollary 3. 
\Box

\vskip3mm

\noindent{\bf Proposition 8.} {\sl 
Fix $\vb\in\Vb$, $B\in\Mb$ and $A=E_{_{h,h+1}}({\bf v},a)$.
Put $l = \max\{i\,|\, b_{h+1, i}\neq 0\}$ and suppose, in addition, 
that $b_{h+1,l}\geq a$ and $\sum_jb_{ij}=\sum_ja_{ji}$.\hfill\break
(a) Then $A\circ B = B + a\,(E_{hl} - E_{h+1, l})=T_{13}$ where
$T=(t_{ijk})$ is the 3-array in $\Tb(A,B)$ such that
$$\left\{\matrix{ 
t_{iik}=b_{ik}-\delta_{i\,h+1}\delta_{k\,l}a,\hfill&\cr
t_{h,h+1,k}=\delta_{k\,l}a,\hfill&\cr
t_{ijk}=0\quad\hfill&{\sl otherwise.}\hfill
}\right.$$
(b) If $b_{hl}=0$ then the projection $p_{_{13,T}}$ is an isomorphism
$O_{_{T}}\buildrel\sim\over\to O_{_{A\circ B}}$. Moreover,

\vskip1mm

\item{\bf .} 
$S_{A\circ B}\subset S_B$ is the subgroup formed by the $w$'s which preserve 
the sets $[A\circ B]_{_{h+1,l}}$ and $[A\circ B]_{_{h,l}},$

\item{\bf .}
$w\cdot S_{A\circ B}\cdot w^{-1}\subset S_A$
if $w\in S_d$ is such that for all $i\neq h$ the following holds
$$w([A]_{_{hh}})=\bigcup_{j\neq l}[A\circ B]_{_{hj}},\quad
w([A]_{_{h\,h+1}})=[A\circ B]_{_{hl}},\quad
w([A]_{_{ij}})=\bigcup_j[A\circ B]_{_{ij}}.$$

\noindent Then, $\Rb^{^{(A)}}_{_\AA}$ and $\Rb^{^{(B)}}_{_\AA}$ 
may be viewed as subrings of $\Rb^{^{(A\circ B)}}_{_\AA}$ and, if
$f\in\Rb^{^{(A)}}_{_\AA}$, $g\in\Rb^{^{(B)}}_{_\AA},$ then
$$f\star g=f\,g\in\Rb^{^{(A\circ B)}}_{_\AA}.$$

\vskip1mm

\noindent (c) If $A - a_{_{h, h-1}}\cdot E_{_{h, h-1}}$ is diagonal,
$l = \min\{ i\,|\, b_{h-1, i}\neq 0\}$ and $b_{h-1,l}\geq a_{h,h-1}$,
then we have similar formulas with 
$A\circ B = B + a_{h, h-1}(E_{hl}- E_{h-1, l})$.
}

\vskip3mm

\proclaim Lemma 9. Let
$A=E_{_{h,h+1}}({\bf v},a)$. For any matrix
$B$ such that $\sum_jb_{ij}=\sum_ja_{ji}$, 
the set ${\bf T}(A,B)$ is in bijective correspondence with
the following set of $n$-tuples
$${\bf S}(A,B)=\{s=(s_1,s_2,...,s_n)\,
|\, 0\leq s_k\leq b_{h+1,k},\quad \sum_ks_k=a\}$$
To an $n$-tuple $s\in{\bf S}(A,B)$ one assigns the following array
$T(s)=(t_{ijk})$ 
$$\left\{\matrix{ 
t_{iik}=b_{ik}\quad\hfill&{\sl if}\quad i\neq h+1,\hfill\cr 
t_{h+1,h+1,k}=b_{h+1, k}-s_k,\hfill&\cr 
t_{h,h+1,k}=s_k,\hfill&\cr
t_{ijk}=0\quad\hfill&{\sl otherwise.}\hfill
}\right.$$
Moreover $T(s)_{13}=B+\sum_j s_j\cdot(E_{hj}-E_{h+1,j})$.

\vskip3mm

\noindent {\sl Proof of the lemma.} An array $T=(t_{ijk})$
belongs to ${\bf T}(A,B)$ if and only if
$$\sum_kt_{ijk}=a_{ij}\quad{\rm and}\quad\sum_it_{ijk}=b_{jk}.$$
In particular,
$$\left\{\matrix{ 
t_{ijk}=0\quad\hfill&{\rm if}\quad (i,j)\neq (i,i),(h,h+1),\hfill\cr
t_{jjk}=b_{jk}\quad\hfill&{\rm if}\quad j\neq h+1,\hfill\cr
t_{h, h+1, k}+t_{h+1,h+1,k}=b_{h+1,k}.\hfill&
}\right.$$
Thus $T\in\Tb(A,B)$ if and only if 
$$\left\{\matrix{ 
t_{ijk}=0\quad\hfill&{\rm if}\quad (i,j)\neq (i,i),(h,h+1),\hfill\cr
t_{jjk}=b_{jk}\quad\hfill&{\rm if}\quad j\neq h+1,\hfill\cr
t_{h+1, h+1, k} = b_{h+1,k}-t_{h,h+1,k},\hfill&\cr
t_{h, h+1, k}\leq b_{h+1,k},\hfill&\cr
\sum_kt_{h, h+1, k}=a_{h,h+1},\hfill&\cr
\sum_kt_{h+1, h+1, k}=a_{h+1,h+1},\hfill&\cr
\sum_kb_{jk}=a_{jj}\quad\hfill&{\rm if}\quad j\neq h+1.\hfill
}\right.$$
The last two equations follow directly from $\sum_jb_{ij}=\sum_ja_{ji}$.
The conditions above mean that there exists an $n$-tuple 
$s\in {\bf S}(A,B)$ such that
$$\left\{\matrix{ 
t_{iik}=b_{ik}\quad\hfill&{\rm if}\quad i\neq h+1,\hfill\cr
t_{h+1,h+1,k}=b_{h+1, k}-s_k,\hfill\cr
t_{h,h+1,k}=s_k,\hfill\cr
t_{ijk}=0\quad\hfill&{\rm otherwise.}\hfill
}\right.$$
\Box

\vskip.3cm

\noindent {\sl Proof of the proposition.} For $s\in{\bf S}(A,B)$ we put
$C(s)= T(s)_{13}$. Observe that, if
$s,s'\in{\bf S}(A,B)$ then we have (see Section 1.2)
$$C(s')\preceq C(s)\Leftrightarrow\quad
\left\{\matrix{
\sum_{r\geq j}s'_r\leq\sum_{r\geq j}s_r\quad\hfill&\forall j>h,\hfill\cr\cr
\sum_{r\leq j}s'_r\geq\sum_{r\leq j}s_r\quad\hfill&\forall j\leq h.
}\right.$$
Thus, since $\sum_ks'_k=\sum_ks_k(=a)$, 
$$C(s')\preceq C(s)\Leftrightarrow\quad\forall j,\quad
\sum_{r\geq j}s'_r\leq\sum_{r\geq j}s_r.$$
Observe on the other hand that
$b_{h+1,k}=0$ if $k>l$. Hence $s_k=0$ if $k>l$, for any
$s\in{\bf S}(A,B)$. Since $b_{h+1,l}\geq a$, the set 
$\{C(s)\,|\,s\in{\bf S}(A,B)\}$ has a greatest
element with respect to $\preceq$ :
it is labelled by the $n$-tuple
$s=a\cdot{\bf e}_l\in{\bf S}(A,B)$.
Since the map ${\bf m}$ is monotone (see Section 1.2), 
it follows that $A\circ B=C(a\cdot\eb_l)$.
Then, the formula in Lemma 9 yields 
$$A\circ B=B+a(E_{hl}-E_{h+1,l})$$
and the corresponding array is the array $T$
defined in Proposition 8 (a). We obtain 
$$O_{_T}=\{(D,D',D'')\in p_{_{23}}^{^{-1}}O_{_B}\,|
\,D'_i\subseteq D_i\quad{\rm and}\quad\dim(D_i\cap D''_j)=\dim(D'_i\cap D''_j)+
\delta_{_{hi}}\delta_{_{j\geq l}}a\}.$$
For any $(D,D'')\in O_{_{A\circ B}}$, 
$$p_{_{13, T}}^{-1}(D,D'')=\{V\,|\,D_{h-1}\i V\i D_h\quad{\rm and}\quad 
{\rm dim}(D_h\cap D''_j)={\rm dim}(V\cap D''_j)+a\delta_{j\geq l}\}.$$
Thus, the fiber is isomorphic to 
the Grassmannian of codimension $a$ subspaces in
the $(a+b_{hl})$-dimensional vector space
$${D_h\cap D''_l\over(D_{h-1}\cap D''_l)+(D_h\cap D''_{l-1})}.$$
In particular, if $b_{hl}=0$ the projection $p_{_{13,T}}$ induces an isomorphism
$O_{_{T}}\buildrel\sim\over\to O_{_{A\circ B}}$.
The maps
$$O_{_{A\circ B}}\buildrel\sim\over\to O_{_{T}}\longrightarrow O_{_A}
\quad{\rm and}\quad
O_{_{A\circ B}}\buildrel\sim\over\to O_{_{T}}\longrightarrow O_{_B}$$
are $GL_d$-equivariant maps between $GL_d$-orbits.
The corresponding pull-back morphisms in equivariant $K$-theory coincide,
via the induction, with the restriction of representations of the 
isotropy subgroups. Then apply Corollary 3.
\Box 

\vskip.3cm

\noindent {\bf 1.5. Generators of the convolution algebra.} 
We use the notations introduced in the previous sections.
Recall that if $A\in{\bf M}$ is diagonal or is of type 
$E_{_{i\,i\pm 1}}({\bf v},a)$ then the $GL_d$-orbit $O_{_A}$ is closed in $F^2$.

\vskip3mm

\noindent {\bf Proposition 10.} {\sl The convolution $\KK$-algebra $(K^G(Z)_\KK,\star)$
is generated by the classes of sheaves supported on the irreducible components
$\overline Z_{_A}$ such that $A\in{\bf M}$ is a diagonal matrix or 
a matrix of type $E_{_{i\,i\pm 1}}({\bf v},1)$.}

\vskip3mm

\noindent {\sl Proof.} Fix a matrix $C\in{\bf M}(\vb,\wb)$.
The proof goes by induction on
$$l(C)=\sum_{i\neq j}\pmatrix{|i-j|+1\cr 2} c_{ij}.$$
If $l(C)=0$ or $1$ the matrix $C$ is of the prescribed type. 
Suppose that $l(C)>1$. 
Then $C$ is not diagonal. Suppose for
instance that $C$ is not a lower-triangular matrix. Put 
$$(h,l)={\rm max}\{(i,j)\,|\,1\leq i<j\leq n,\quad c_{ij}\neq 0\},$$
with respect to the right-lexicographic order (see Section 1.1). Put 
$$B=C+c_{_{hl}}(E_{_{h+1,l}}-E_{_{hl}})\quad{\rm and}\quad
A=E_{_{h,h+1}}(\vb-c_{_{hl}}\cdot{\bf e}_{_h},c_{_{hl}}).$$
Then $l={\rm max}\{i\,|\,b_{h+1,i}\neq 0\}$ and 
$b_{h+1,l}=c_{hl}+c_{h+1,l}\geq c_{hl}=a_{h,h+1}.$
Moreover,
$$\sum_ja_{ji}=\sum_jc_{ij}-\delta_{ih}c_{hl}+\delta_{i,h+1}c_{hl}
=\sum_jb_{ij}$$
and $b_{hl}=0$.
Thus we are in the situation of Proposition 8 (b).
In particular, $A\circ B=C$.
Given $f\in\Rb^{^{(A)}}_{_\AA}$ and $g\in\Rb^{^{(B)}}_{_\AA}$ we get
$$f\star g=f\,g\in\Rb^{^{(C)}}_{_\AA},$$
modulo the identifications described in Proposition 8 (b).
Thus the surjectivity of the map 
$\star\,:\,\Rb^{^{(A)}}_{_\AA}\otimes\Rb^{^{(B)}}_{_\AA}
\longrightarrow\Rb^{^{(C)}}_{_\AA}$
follows from the surjectivity of the map
$$\Rb^{(a_{h\,h+1})}_{_\AA}\otimes\Rb^{(b_{h+1\,l})}_{_\AA}\longrightarrow
\Rb^{(c_{h\,l})}_{_\AA}\otimes\Rb^{(c_{h+1\,l})}_{_\AA},\qquad
f\otimes g\mapsto i(f)\,j(g),$$
where $i$ is the obvious isomorphism 
$\Rb^{(a_{h,h+1})}_{_\AA}\simeq\Rb_{_\AA}^{(c_{hl})}$ 
and $j$ is the inclusion 
$\Rb^{(b_{h+1,l})}_{_\AA}\to\Rb^{(c_{hl})}_{_\AA}\otimes
\Rb^{(c_{h+1\,l})}_{_\AA}$ 
(for any partition $[d]=I_1\cup I_2$ the algebra $\Rb^{^{(I_1,I_2)}}$ 
is generated by $\Rb^{^{(d)}}$ and $\Rb^{^{(I_1)}}$).
On the other hand $l(B)<l(C)$.
To complete the proof it suffices to remark that any sheaf supported on
$Z_{E_{i\,i\pm 1}(\vb,a)}$ can be obtain by the convolution
product of $a$ sheaves supported respectively on 
$$Z_{E_{i\,i\pm 1}(\vb+k\eb_{i\pm 1},1)},\qquad k=0,1,...,a-1,$$
using the computations done in example 1.4.
\Box

\vskip1cm

\noindent{\bf 2. The polynomial representation and the convolution algebra.}

\vskip 3mm

\noindent Fix $n,d\geq 1$.

\vskip 3mm

\noindent{\bf 2.1. The Drinfeld new presentation.} 
The affine quantum group of type ${\got{gl}}\,(n)$ 
with trivial central charge
is a unital associative algebra over 
$\KK=\CC(q)$ 
with generators
$$\Eb_{i,k},\quad\Fb_{i,k},\quad\Kb_{j,l},\quad\Kb_{j}^{^{\pm 1}},
\quad i\in [n-1],\quad j\in [n],\quad k,l\in\ZZ,\quad l\neq 0.$$
The relations are expressed in terms of the formal series
$$\Eb_i(z)=\sum_{k\in\ZZ}{\Eb_{i,k}\cdot z^{^{-k}}},\quad
\Fb_i(z)=\sum_{k\in\ZZ}{\Fb_{i,k}\cdot z^{^{-k}}},\quad
\Kb_i^{^\pm}(z)=\Kb^{^{\pm 1}}_{i}+
\sum_{l>0}{\Kb_{i,\pm l}\cdot z^{^{\mp l}}},$$
as follows
$$\Kb_{i}^{^{-1}}\cdot\Kb_{i}=\Kb_{i}\cdot\Kb_{i}^{^{-1}}=1,\quad 
[\Kb^{^\pm}_i(z),\Kb^{^\pm}_j(w)]=[\Kb^{^\pm}_i(z),\Kb^{^{\mp}}_i(w)]=0,
\leqno{\rm (a)}$$
$$\theta_1(q^{^{-1}} z/w)\cdot\Kb^{^\pm}_i(z)\cdot
\Kb^{^\mp}_j(w)=\theta_1(q^{^{-1}} z/w)\cdot
\Kb^{^\mp}_j(w)\cdot\Kb^{^\pm}_i(z),\quad{\rm if}\quad i<j\leqno{\rm (b)}$$
$$\Kb^{^\pm}_j(z)\cdot\Eb_i(w) =\theta_{c_{ij}}
(q^{^{c_{ij}}} z/w)\cdot\Eb_j(w)\cdot\Kb^{^\pm}_i(z),
\leqno{\rm (c)}$$
$$\Kb^{^\pm}_j(z)\cdot\Fb_i(w) =\theta_{-c_{ij}}
(q^{^{c_{ij}}} z/w)\cdot\Fb_i(w)\cdot\Kb^{^\pm}_j(z),
\leqno{\rm (d)}$$
$$[\Eb_i(z),\Fb_j(w)]=(q-q^{^{-1}})\cdot\delta_{ij}\cdot
\leqno{\rm (e)}$$
\hfill$\biggl(\delta(z/w)\cdot\Kb^{^+}_{i+1}(w)/
\Kb^{^+}_{i}(w)-
\delta(z/w)\cdot\Kb^{^-}_{i+1}(z)/\Kb^{^-}_i(z)\biggr),$
\hfill\break
$$\Eb_i(z)\cdot\Eb_j(w)=\theta_{m_{ij}}
(q^{^{i-j}}z/w)\cdot\Eb_j(w)\cdot\Eb_i(z),\leqno{\rm (f)}$$
$$\Fb_i(z)\cdot\Fb_j(w)=\theta_{-m_{ij}}
(q^{^{i-j}}z/w)\cdot\Fb_j(w)\cdot\Fb_i(z),\leqno{\rm (g)}$$
$$\{(\Eb_i(z_{_1})\cdot\Eb_i(z_{_2})\cdot\Eb_{j}(w)-
(q+q^{^{-1}})\cdot\Eb_i(z_{_1})\cdot\Eb_{j}(w)\cdot\Eb_i(z_{_{2}})+
\Eb_{j}(w)\cdot\Eb_i(z_{_{1}})\cdot\Eb_i(z_{_2})\}+\leqno{\rm (h)}$$
$$+\{z_{_1}\leftrightarrow z_{_2}\}=0,
\qquad{\rm if}\quad |i-j|=1,$$
$$\{\Fb_i(z_{_1})\cdot\Fb_i(z_{_2})\cdot\Fb_{j}(w)-
(q+q^{^{-1}})\cdot\Fb_i(z_{_1})\cdot\Fb_{j}(w)\cdot\Fb_i(z_{_{2}})+
\Fb_{j}(w)\cdot\Fb_i(z_{_{1}})\cdot\Fb_i(z_{_2})\}+\leqno{\rm (i)}$$
$$+\{z_{_1}\leftrightarrow z_{_2}\}=0,
\qquad{\rm if}\quad |i-j|=1,$$
$$[\Fb_i(z),\Fb_j(w)]=[\Eb_i(z),\Eb_j(w)]=0,
\qquad{\rm if}\quad |i-j|>1,\leqno{\rm (j)}$$
where $\delta(z) = \sum_{_{n=-\infty}}^{^\infty}z^n$,
$\theta_m(z)={q^{^m}\cdot z -1\over z-q^{^m}}$ and
$c_{ij}$, $m_{ij}$, are the entries of the following $n\times n$-matrices 
\hfill\break

\hfill$C=\pmatrix{
-1&1&0&&0\cr
0&-1&1&&0\cr
0&0&-1&\cdots&0\cr
&&\vdots&\ddots&\vdots\cr
0&0&0&\cdots&-1\cr}
\qquad{\rm and}\qquad
M=-C-C^t.$
\hfill\break

\noindent See [D] and [DF] for more details. In particular 
the relations (c), (d), (f) and (g) are written formally, and in practice
one should first multiply both sides by the denominator  
of the function $\theta_m$ since formal currents may have zero divisors.

\vskip3mm

\noindent{\bf 2.2. The polynomial representation.} 
We keep the same notations as in the first chapter. 
Thus $\Rb=\CC[x_1^{\pm 1},x_2^{\pm 1},...,x_d^{\pm 1}]$
and for any composition 
$\vb\in\Vb$ the algebra $\Rb^{(\vb)}$ is identified with the subalgebra
of $S_{\vb}$-invariant polynomials in $\Rb$.
Set $\Kb=\bigoplus_{\vb\in\Vb}\Rb^{(\vb)}$.
For any subset $I\i[d]$ set
$$\Theta_{_I}(z)=\prod_{m\in I}\theta_1(z/x_m).$$
Recall that if $\vb\in\Vb$ and $i\in[n]$ then 
$\overline v_i=v_1+v_2+...+v_i$ and
$[\vb]_i=[1+\overline v_{i-1},\overline v_i]$. 
Let 
$$\hat\Eb_i(z),\hat\Fb_i(z)\in\bigoplus_{\vb\in\Vb}
{\rm Hom}(\Rb^{^{(\vb)}}_{_\AA},\Rb_{_\AA}^{^{(\vb\pm\eb_i\mp\eb_{i+1})}})
[[z^{^{\pm 1}}]]$$
be the following operators
$$\hat\Eb_i(z)(f)=\cases{
{\displaystyle
(q-q^{^{-1}})\,{\got S}_\vb^{\vb+\eb_i-\eb_{i+1}}
\biggl(f\cdot\delta(x_{1+\overline v_i}/z)\cdot
\Theta_{_{[\vb]_i}}(qx_{1+\overline v_i})
\biggr)}\quad &if $\vb+\eb_i-\eb_{i+1}\in\Vb$,\cr
0& else,
}$$
$$\hat\Fb_i(z)(f)=\cases{
{\displaystyle
(q-q^{^{-1}})\,{\got S}_\vb^{\vb-\eb_i+\eb_{i+1}}
\biggl(f\cdot\delta(x_{\overline v_i}/z)\cdot
\Theta_{_{[\vb]_{i+1}}}(x_{\overline v_i}/q)^{-1}
\biggr)}\quad &if $\vb-\eb_i+\eb_{i+1}\in\Vb$,\cr
0& else.
}$$
Moreover let $\hat\Kb^\pm_i(z)\in{\rm End}(\Kb_{_\AA})[[z^{^{\pm 1}}]]$ be 
the operator whose Fourier coefficients act on $\Rb^{^{(\vb)}}$ 
by multiplication by the corresponding coefficient of the expansion
at $z=\infty, 0$ respectively of the rational function
$$\hat\Kb_{i,\vb}(z)= \Theta_{_{[1,\overline v_{i-1}]}}(qz)\cdot
\Theta_{_{[1+\overline v_i,d]}}(z/q).$$
This rational function differs from the one used in [GV]
since we consider here type $\widehat{\got gl}(n)$
instead of $\widehat{\got sl}(n)$.
Denote by $\hat\Eb_{i,k}$, $\hat\Fb_{i,k}$,
$\hat\Kb_{i,k}$ and $\hat\Kb_i^{\pm 1}$ if $k=0$
the Fourier coefficient of $z^{-k}$ in the series 
$\hat\Eb_i(z)$, $\hat\Fb_i(z)$,
$\hat\Kb_i^\pm(z)$ (see Section 2.1).
The main result of this section is the following theorem.

\vskip3mm

\noindent{\bf Proposition 11.} {\sl The map
$$\Eb_i(z)\mapsto\hat\Eb_i(z),\qquad
\Fb_i(z)\mapsto\hat\Fb_i(z),\qquad
\Kb_i^\pm(z)\mapsto\hat\Kb_i^\pm(z),
$$
extends uniquely to a representation of $\Ub$ on $\Kb_{_\KK}$.}

\vskip3mm

\noindent{\sl Proof.} Our proof is an extension to the affine case of
the computations done in [V]. We will prove that the operators $\hat\Eb_i$,
$\hat\Fb_i$ and $\hat\Kb_i$ satisfy the relations given in Section 2.1.
The first and the second one are immediate. 
To verify the others we need some more notation.
First, given a partition $I=(I_1,I_2,...,I_n)$ of the set 
$[d]$ into $n$ subsets, for any $x\in(\CC^\times)^d$
let $x_I\in\CC^{\times (I)}$ be its projection.
Second, to any $k\in [d]$ we associate two operators,
$\tau_k^\pm$, on the set of partitions of $[d]$ in such a way that
$\tau_k^\pm I$ is the partition with $k$ shifted from one piece of $I$
to the next (resp. the previous) one.
Then, for any $\vb\in\Vb$ such that $\vb\mp\eb_i\pm\eb_{i+1}\in\Vb$
and any $f\in\Rb^{^{(\vb\mp\eb_i\pm\eb_{i+1})}}_{_\AA}$ we get respectively
$$(\hat\Eb_i(z)f)(x_{[\vb]})= (q-q^{^{-1}})\, \sum_{k\in[\vb]_i}
f(x_{\tau^+_k[\vb]})
\delta(x_k/z)\Theta_{_{[\vb]_i\setminus\{k\}}}(qx_k),$$
$$(\hat\Fb_i(z)f)(x_{[\vb]})= (q-q^{^{-1}})\, \sum_{k\in[\vb]_{i+1}}
f(x_{\tau_k^-[\vb]})
\delta(x_k/z)\Theta_{_{[\vb]_{i+1}\setminus\{k\}}}(x_k/q)^{-1}.
$$

\noindent{\bf .} Let us prove (e). First observe that
$$\hat\Kb_{i+1,\vb}(z)\hat\Kb_{i,\vb}(z)^{-1}=
\Theta_{_{[\vb]_i}}(qz)\Theta_{_{[\vb]_{i+1}}}(z/q)^{-1}.$$
The values of
$(q-q^{^{-1}})^{^{-2}}\hat\Eb_i(z)\hat\Fb_i(w)f$ and 
$(q-q^{^{-1}})^{^{-2}}\hat\Fb_i(w)\hat\Eb_i(z)f$ 
at $x_{_{[\vb]}}$ are respectively
$$\sum_{k\in[\vb]_i\atop l\in[\vb]_{i+1}\cup\{k\}}
f(x_{\tau^-_l\tau^+_k[\vb]})\,\delta(x_k/z)\,\delta(x_l/w)\,
\Theta_{_{[\vb]_i\setminus\{k\}}}(qx_k)
\Theta_{_{[\vb]_{i+1}\cup\{k\}\setminus\{l\}}}(x_l/q)^{-1},$$
$$\sum_{k\in[\vb]_i\cup\{l\}\atop l\in[\vb]_{i+1}}
f(x_{\tau^+_k\tau^-_l[\vb]})\,\delta(x_k/z)\,\delta(x_l/w)\,
\Theta_{_{[\vb]_i\cup\{l\}\setminus\{k\}}}(qx_k)
\Theta_{_{[\vb]_{i+1}\setminus\{l\}}}(x_l/q)^{-1}.$$
Their difference is equal to :
$$f(x_{[\vb]})\,
\sum_{k\in[\vb]_i}\delta(x_k/z)\,\delta(x_k/w)\,
\Theta_{_{[\vb]_i\setminus\{k\}}}(qx_k)\Theta_{_{[\vb]_{i+1}}}(x_k/q)^{-1}-$$
$$-f(x_{[\vb]})\,
\sum_{k\in[\vb]_{i+1}}\delta(x_k/z)\,\delta(x_k/w)\,
\Theta_{_{[\vb]_i}}(qx_k)\Theta_{_{[\vb]_{i+1}\setminus\{k\}}}(x_k/q)^{-1}.$$
Set $A(x)=\prod_{k\in[\vb]_i\cup[\vb]_{i+1}}(x-x_k)$ and
$B(x)=\prod_{k\in[\vb]_i}(qx-q^{^{-1}}x_k)
\prod_{k\in[\vb]_{i+1}}(q^{^{-1}}x-qx_k)$.
Then
$$([\hat\Eb_i(z),\hat\Fb_i(w)]f)(x_{_{[\vb]}})=(q-q^{^{-1}})^{^2}
\sum_{k\in[\vb]_i\cup[\vb]_{i+1}}
{\delta(x_k/z)\,\delta(x_k/w)\,x_k^{^{-1}}\,B(x_k)\over A'(x_k)}
f(x_{_{[\vb]}}).$$
The right hand side in (e) is then obtained by a computation of residues.
If $i\neq j\pm 1$ the preceding formulas for $\hat\Eb_i(z)$ and
$\hat\Fb_j(w)$ give immediately $[\hat\Eb_i(z),\hat\Fb_j(w)]=0$.
Similarly the values of $(q-q^{^{-1}})^{^{-2}}\hat\Eb_i(z)\hat\Fb_{i+1}(w)f$ and 
$(q-q^{^{-1}})^{^{-2}}\hat\Fb_{i+1}(w)\hat\Eb_i(z)f$ at $x_{_{[\vb]}}$ 
are both equal to
$$\sum_{k\in[\vb]_i\atop l\in[\vb]_{i+2}}
f(x_{\tau_k^+\tau_l^-[\vb]})\,\delta(x_k/z)\,\delta(x_l/w)
\Theta_{_{[\vb]_i\setminus\{k\}}}(qx_k)
\Theta_{_{[\vb]_{i+2}\setminus\{l\}}}(x_l/q)^{-1},$$
whereas the values of $(q-q^{^{-1}})^{^{-2}}\hat\Eb_i(z)\hat\Fb_{i-1}(w)f$ and 
$(q-q^{^{-1}})^{^{-2}}\hat\Fb_{i-1}(w)\hat\Eb_i(z)f$ at $x_{_{[\vb]}}$ are 
both equal to
$$\sum_{k,l\in[\vb]_i\atop k\neq l}
f(x_{\tau_k^+\tau_l^-[\vb]})\,\delta(x_k/z)\,\delta(x_l/w)\,
\Theta_{_{[\vb]_i\setminus\{k\}}}(qx_k)
\Theta_{_{[\vb]_i\setminus\{k,l\}}}(x_l/q)^{-1}.$$
Thus $[\hat\Eb_i(z),\hat\Fb_j(w)]=0$ whenever $i\neq j$.

\vskip1mm

\noindent{\bf .} Let us prove (f). The
value of $(q-q^{^{-1}})^{^{-2}}\hat\Eb_k(z)\hat\Eb_l(w)f$ at 
$x_{_{[\vb]}}$ in the three cases $(k,l)=(i,i),(i,i+1),(i+1,i)$
is respectively
$$\sum_{k,l\in[\vb]_i\atop k\neq l}
f(x_{\tau_k^+\tau_l^+[\vb]})\,\delta(x_k/z)\,\delta(x_l/w)\,
\Theta_{_{[\vb]_i\setminus\{k\}}}(qx_k)\,
\Theta_{_{[\vb]_i\setminus\{k,l\}}}(qx_l),$$
$$\sum_{k\in[\vb]_i\atop l\in[\vb]_{i+1}\cup\{k\}}
f(x_{\tau_k^+\tau_l^+[\vb]})\,\delta(x_k/z)\,\delta(x_l/w)\,
\Theta_{_{[\vb]_i\setminus\{k\}}}(qx_k)\,
\Theta_{_{[\vb]_{i+1}\cup\{k\}\setminus\{l\}}}(qx_l),$$
$$\sum_{k\in[\vb]_{i+1}\atop l\in[\vb]_i}
f(x_{\tau_k^+\tau_l^+[\vb]})\,\delta(x_k/z)\,\delta(x_l/w)\,
\Theta_{_{[\vb]_{i+1}\setminus\{k\}}}(qx_k)\,
\Theta_{_{[\vb]_i\setminus\{l\}}}(qx_l).$$
The first equality gives
$$(q^{^{-1}}z-qw)\hat\Eb_i(z)\hat\Eb_i(w)=
(qz-q^{^{-1}}w)\hat\Eb_i(w)\hat\Eb_i(z),$$
the other two give 
$$(w-z)\hat\Eb_i(z)\hat\Eb_{i+1}(w)=
(qw-q^{^{-1}}z)\hat\Eb_{i+1}(w)\hat\Eb_i(z)$$
(note that we have used the basic formula $(z-w)\,\delta(z)\,\delta(w)=0$).
When $|i-j|>1$ the relation (f) is immediate.
The proof of (g) is similar.

\vskip1mm

\noindent{\bf .} Let us prove (c). The 
value of $(q-q^{^{-1}})^{^{-2}}\hat\Eb_k(z)\hat\Eb_l(w)f$ at 
$x_{_{[\vb]}}$ in the three cases $(k,l)=(i,i),(i,i+1),(i+1,i)$
is respectively given by
$$\sum_{k,l\in[\vb]_i\atop k\neq l}
f(x_{\tau_k^+\tau_l^+[\vb]})\,\delta(x_k/z)\,\delta(x_l/w)\,
\Theta_{_{[\vb]_i\setminus\{k\}}}(qx_k)\,
\Theta_{_{[\vb]_i\setminus\{k,l\}}}(qx_l),$$

\vskip5mm

\noindent$(\hat\Eb_i(z)\hat\Kb_i(w)f)(x_{_{[\vb]}})=$\hfill\break
$$=(q-q^{^{-1}})\, \sum_{k\in[\vb]_i}
f(x_{\tau_k^+[\vb]})
\delta(x_k/z)\Theta_{_{[\vb]_i\setminus\{k\}}}(qx_k)
\Theta_{_{[1,\overline v_{i-1}]}}(qw)\Theta_{_{[1+\overline v_i,d]}}(w/q)
\theta_{_1}(q^{^{-1}}w/x_k),$$

\noindent and

\vskip5mm

\noindent$(\hat\Kb_i(w)\hat\Eb_i(z)f)(x_{_{[\vb]}})=$\hfill\break
$$=(q-q^{^{-1}})\, \sum_{k\in[\vb]_i}
f(x_{\tau_k^+[\vb]})
\delta(x_k/z)\Theta_{_{[\vb]_i\setminus\{k\}}}(qx_k)
\Theta_{_{[1,\overline v_{i-1}]}}(qw)\Theta_{_{[1+\overline v_i,d]}}(w/q)
.$$

\noindent Since $\delta(x_k/z)\theta_{_1}(q^{^{-1}}w/x_k)=
\delta(x_k/z)\theta_{_1}(q^{^{-1}}w/z)$, we thus obtain 

$$\hat\Eb_i(z)\hat\Kb_i(w)=
\theta_{_1}(q^{^{-1}}w/z)\hat\Kb_i(w)\hat\Eb_i(z).$$

\noindent Similarly we have

\vskip5mm

\noindent$(\hat\Eb_i(z)\hat\Kb_{i+1}(w)f)(x_{_{[\vb]}})=$\hfill\break
$$=(q-q^{^{-1}})\,\sum_{k\in[\vb]_i}
f(x_{\tau_k^+[\vb]})
\delta(x_k/z)\Theta_{_{[\vb]_i\setminus\{k\}}}(qx_k)
\Theta_{_{[1,\overline v_{i}]}}(qw)\Theta_{_{[1+\overline v_{i+1},d]}}(w/q)
\theta_{_{-1}}(qw/x_k),$$

\noindent and

\vskip5mm

\noindent$(\hat\Kb_{i+1}(w)\hat\Eb_i(z)f)(x_{_{[\vb]}})=$\hfill\break
$$=(q-q^{^{-1}})\,\sum_{k\in[\vb]_i}
f(x_{\tau_k^+[\vb]})
\delta(x_k/z)\Theta_{_{[\vb]_i\setminus\{k\}}}(qx_k)
\Theta_{_{[1,\overline v_i]}}(qw)\Theta_{_{[1+\overline v_{i+1},d]}}(w/q)
.$$

\noindent Thus,

$$\hat\Eb_i(z)\hat\Kb_{i+1}(w)=
\theta_{_{-1}}(qw/z)\hat\Kb_{i+1}(w)\hat\Eb_i(z).$$

\noindent The other cases are obvious and (d) is proved in the same way. 

\vskip1mm

\noindent{\bf .} Let us prove formula (j). If $|i-j|>2$ then
$(q-q^{^{-1}})^{^{-2}}\hat\Eb_i(z)\hat\Eb_j(w)f(x_{[\vb]})$ equals 
$$\sum_{k\in[\vb]_i\atop l\in[\vb]_j}
f(x_{\tau_l^+\tau_k^+[\vb]})\,\delta(x_k/z)\,\delta(x_l/w)\,
\Theta_{_{[\vb]_i\setminus\{k\}}}(qx_k)\,
\Theta_{_{[\vb]_j\setminus\{l\}}}(qx_l).$$
Switching the factors $\hat\Eb_i(z)$ and $\hat\Eb_j(w)$ yields the same
expression.

\vskip1mm

\noindent{\bf .} Let us prove the remaining Serre relations.
We will only prove the formula (h) for $j=i+1$ since the other
cases are similar. First, let notice that
$(q-q^{^{-1}})^{^3}
\bigr(\hat\Eb_i(z_{_1})\hat\Eb_i(z_{_2})\hat\Eb_{i+1}(w)\,f\bigl)(x_{[\vb]})$,
$(q-q^{^{-1}})^{^3}
\bigr(\hat\Eb_i(z_{_1})\hat\Eb_{i+1}(w)\hat\Eb_i(z_{_2})\,f\bigl)(x_{[\vb]})$ and
$(q-q^{^{-1}})^{^3}
\bigr(\hat\Eb_{i+1}(w)\hat\Eb_i(z_{_1})\hat\Eb_i(z_{_1})\,f\bigl)(x_{[\vb]})$
are equal respectively to
$$\sum_{k\neq l\in[\vb]_i\atop h\in[\vb]_{i+1}\cup\{k,l\}}
A_{k,l,h}(z_1,z_2,w)\,
\Theta_{_{[\vb]_i\setminus\{k\}}}(qx_k)\,
\Theta_{_{[\vb]_i\setminus\{k,l\}}}(qx_l)\,
\Theta_{_{[\vb]_{i+1}\cup\{k,l\}\setminus\{h\}}}(qx_h),$$
$$\sum_{k\neq l\in[\vb]_i\atop h\in[\vb]_{i+1}\cup\{k\}}
A_{k,l,h}(z_1,z_2,w)\,
\Theta_{_{[\vb]_i\setminus\{k\}}}(qx_k)\,
\Theta_{_{[\vb]_{i+1}\cup\{k\}\setminus\{h\}}}(qx_h)\,
\Theta_{_{[\vb]_i\setminus\{k,l\}}}(qx_l),$$
$$\sum_{k\neq l\in[\vb]_i\atop h\in[\vb]_{i+1}}
A_{k,l,h}(z_1,z_2,w)\,
\Theta_{_{[\vb]_{i+1}\setminus\{h\}}}(qx_h)\,
\Theta_{_{[\vb]_i\setminus\{k\}}}(qx_k)\,
\Theta_{_{[\vb]_i\setminus\{l,k\}}}(qx_l),$$
where 
$A_{k,l,h}(y,z,w)=f(x_{\tau_k^+\tau_l^+\tau_h^+[\vb]})\,\delta(x_k/y)\,
\delta(x_l/z)\,\delta(x_h/w).$
Set $I=[\vb]_{i+1}$. 
We are thus reduced to prove that, given $k\neq l\in[\vb]_i$, we have
$$2\sum_{h\in I\cup\{k,l\}}\Theta_{_{I\cup\{k,l\}\setminus\{h\}}}(qx_h)+
2\sum_{h\in I}\Theta_{_{I\setminus\{h\}}}(qx_h)=$$
$$=(q+q^{-1})\sum_{h\in I\cup\{k\}}\Theta_{_{I\cup\{k\}\setminus\{h\}}}(qx_h)+
(q+q^{-1})\sum_{h\in I\cup\{l\}}\Theta_{_{I\cup\{l\}\setminus\{h\}}}(qx_h).$$
For any subset $J\subset [d]$ with $m$ elements, we have
$$\sum_{h\in J}\Theta_{J\setminus\{h\}}(qx_h)=[m]$$
(see [M] for instance).
Thus the equality is a consequence of the (obvious) equality
$$[m+1]+[m-1]=[2]\,[m],\qquad\forall m\in\NN^\times,$$
which follows from the decomposition of the product of
$m$-dimensional and $2$-dimensional
simple representations of SL(2).
\Box

\vskip3mm

\noindent{\bf 2.3. K-theoretic construction of U.} 
Recall that $F$ is the variety of all $n$-step flags in $\CC^d$
and $Z\subset T^*F\times T^*F$ is the Steinberg variety
(see Sections 1.3 and 1.4). The group $G=GL_d\times\CC^\times$ acts
on $Z$ as in Section 1.4 : $GL_d$ acts diagonally and
$z\in\CC^\times$ acts by multiplication by $z^{-2}$ along the fibers.
Given a partition $\vb$ of $d-1$ into $n$ summands and $i\neq j$ we set 
$$E_{_{ij}}(\vb)=E_{_{ij}}({\bf v},1)={\rm diag}({\bf v})+E_{ij}
\in{\bf M}({\bf v}+{\bf e}_i,{\bf v}+{\bf e}_j).$$
If $j=i\pm 1$ the corresponding $GL_d$-orbits in 
$F_{\vb+\eb_i}\times F_{\vb+\eb_{i\pm 1}}$ 
are closed and are given by 
$$O_{E_{_{i\,i+1}}({\bf v})}=\{(D,D')\in F^2\,|\,D'_k\subseteq D_k\quad
{\rm and}\quad {\rm dim}(D_k/D'_k)=\delta_{ki}\},$$
$$O_{E_{_{i\,i-1}}({\bf v})}=\{(D,D')\in F^2\,|\,D_k\subseteq D'_k\quad
{\rm and}\quad {\rm dim}(D'_k/D_k)=\delta_{k\,i-1}\}.$$
Thus both projections
$$F_{{\bf v}+{\bf e}_i}
\buildrel p_1\over\longleftarrow
O_{E_{i\,i\pm 1}({\bf v})}
\buildrel p_2\over\longrightarrow
F_{{\bf v}+{\bf e}_{i\pm 1}}$$ 
are smooth and proper with fibers isomorphic respectively 
to the projective spaces $\PP^{v_i}$ and $\PP^{v_{i\pm 1}}$. 
Recall that the conormal bundle to $O_{E_{i\,i\pm 1}({\bf v})}$
is denoted by $Z_{E_{i\,i\pm 1}({\bf v})}$. Given $k\in\ZZ$ 
let $\Oc_{p_1}(k)$ and $T^*p_1$ be the relative invertible sheaf
$\Oc(k)$ and the relative cotangent sheaf 
along the fibers of $p_1$ respectively.
Let $\Ec_{i,\vb,k}$ be the pull-back to $Z_{E_{i\,i+1}({\bf v})}$ of
the sheaf ${\rm Det}(T^*p_1)\otimes O_{p_1}(k)$ on $O_{E_{i\,i+1}({\bf v})}$. 
Similarly let $\Fc_{i,\vb,k}$ be the pull-back to $Z_{E_{i+1\,i}({\bf v})}$
of the sheaf ${\rm Det}(T^*p_1)\otimes O_{p_1}(k)$ on $O_{E_{i+1\,i}({\bf v})}$. 
Since the orbits $O_{E_{i\pm 1\,i}({\bf v})}$ are closed, the varieties
$Z_{E_{i\pm 1\,i}({\bf v})}$ are irreducible components
of $Z$. Hence, the sheaves $\Ec_{i,\vb,k}$ and $\Fc_{i,\vb,k}$ 
may be viewed as sheaves on $Z$. They come equiped with the natural
$G$-equivariant structure. We set
$$\Ec_{i,k}=\sum_\vb(-q)^{-v_i}[\Ec_{i,\vb,k}]
\quad{\rm and}\quad
\Fc_{i,k}=\sum_\vb(-q)^{-v_{i+1}}[\Fc_{i,\vb,k}]$$ 
in $K^G(Z)$. Form the corresponding generating functions
$$\Ec_i(z)=\sum_k\Ec_{i,k}z^{^{-k}},\quad\Fc_i(z)=\sum_k\Fc_{i,k}z^{^{-k}}\quad
\in K^G(Z)[[z,z^{^{-1}}]].$$
We have introduced in Section 2.2 some generating series $\hat\Kb_i^\pm(z)$ 
whose Fourier coefficients, $\hat\Kb^\pm_{i,k}$ ($k\in\mp\NN$), are 
homogeneous Laurent polynomials in $x_{_1},x_{_2},...,x_{_d}$.
Since $T^*F$ is isomorphic to the conormal bundle of the diagonal
$\Delta\subset F\times F$, the cotangent bundle $T^*F$
is naturally identified with
an irreducible component of $Z$, say $Z_{_\Delta}$.
Under the map 
$$\Kb_{_\AA}
\buildrel\sim\over\longrightarrow
K^G(F)
\buildrel\sim\over\longrightarrow
K^G(T^*F)
\buildrel\sim\over\longrightarrow
K^G(Z_{_\Delta})
\longrightarrow
K^G(Z)$$
the $\hat\Kb_{i,k}$'s can be viewed
as classes of equivariant sheaves on $Z$ supported by
$Z_{_\Delta}$. Let $\Kc_{i,k}$ denote
these classes and write 
$\Kc_i(z)=\sum_k\Kc_{i,k}z^{^{-k}}\in K^G(Z)[[z,z^{^{-1}}]]$.
The algebra $K^G(Z)$ acts by convolution on $K^G(T^*F)\simeq\Kb_{_\AA}$.
The following result is stated in [GV, Proposition 7.7].

\vskip3mm

\noindent{\bf Lemma 12.} {\sl  The convolution action of $\Ec_{i,k}, \Fc_{i,k},
\Kc_{i,k}\in K^G(Z)$ on $K^G(T^*F)=\Kb_{_\AA}$
is given precisely by the operators 
$\hat\Eb_{i,k}$, $\hat\Fb_{i,k}$ and $\hat\Kb_{i,k}$.
}

\vskip3mm

\noindent{\sl Proof.} See [V, Theorem 10]. The formula in the proof
of Proposition 6 gives
$$[Tp_1]=\sum_{\bar v_{i-1}<t\leq\bar v_i}x_{1+\bar v_i}/x_t.$$
Thus, we obtain
$$\left\{\matrix{
[\Ec_{i,\vb,k}]=x_{1+\bar v_i}^k
\prod_{\bar v_{i-1}<t\leq\bar v_i}x_{1+\bar v_i}^{-1}x_t,\hfill\cr\cr
\sum_i(-1)^i[\Lambda^iTp_1]=
\prod_{\bar v_{i-1}<t\leq\bar v_i}(1-q^2x_{1+\bar v_i}/x_t).
}\right.$$
The result follows from Corollary 4 and Proposition 6.
\Box

\vskip3mm

\noindent The following result is stated in [GV, Lemma 7.5].

\vskip3mm

\noindent{\bf Lemma 13.} {\sl The representation of $K^G(Z)$ on 
$K^G(T^*F)$ by convolution is faithful.}

\vskip3mm

\noindent{\sl Proof.} See [CG, claim 7.6.7].
\Box

\vskip3mm

\noindent The following result is stated in [GV, Theorem 7.9.].

\vskip3mm

\proclaim Theorem 1. The map
$\Eb_{i,k}\mapsto\Ec_{i,k},\quad\Fb_{i,k}\mapsto\Fc_{i,k},
\quad\Kb_{i,k}\mapsto\Kc_{i,k},$
extends uniquely to a morphism of algebras $\Psi\,:\,\Ub\to K^G(Z)_{_\KK}$. 

\vskip3mm

\noindent{\sl Proof.} Immediate since Lemmas 12, 13 and 
Proposition 11 imply that the elements 
$\Ec_{i,k}, \Fc_{i,k}, \Kc_{i,k}$
of the algebra $K^G(Z)$ 
satisfy the relations in Section 2.1. 
\Box

\vskip3mm

\noindent{\bf 2.4. Surjectivity of the map $\Psi$.}
For any $t\in\CC^\times$ let $\Ub_t$ and $K^G(Z)_t$ be the 
$\CC^\times$-algebras obtained by specializating
$\Ub$ and $K^G(Z)$ to $q=t$. The map $\Psi$ restricts to an
algebra homomorphism 
$$\Psi_t\,:\,\Ub_t\rightarrow K^G(Z)_t.$$

\vskip3mm

\noindent{\bf Theorem 2.} {\sl If $t$ is not a root of unity
the map $\Psi_t\,:\,\Ub_t\to K^G(Z)_t$ is surjective.}

\vskip3mm

\noindent 
Fix $t\in\CC^\times$ which is not a root of unity.
For any $\vb\in\Vb$ set $$\Ib_\vb=\prod_{_{i=1}}^{^{n-1}}
\prod_{\buildrel m=0\over{\scriptscriptstyle m\neq d-v_i}}^{^d}
{\Kb_i-t^m\over t^{d-v_i}-t^m}\in\Ub_t.$$

\vskip1mm

\noindent{\bf Lemma 14.} {\sl
$\Psi_t(\Ib_\ub)\star K^G(Z_{\vb,\wb})_t\neq 0$ 
if and only if $\ub=\vb$. Similarly 
$K^G(Z_{\vb,\wb})_t\star\Psi_t(\Ib_\ub)\neq 0$ 
if and only if $\ub=\wb$.
}
\Box

\vskip3mm

\noindent{\bf Lemma 15.} {\sl Fix a partition $\vb\in\Vb$. 
If $t\in\CC^\times$ is not a root of unity, then
the Fourier coefficients of the series
$\hat\Kb_{i,\vb}(z)$'s generate the algebra $\Rb^{(\vb)}$.}

\vskip3mm

\noindent{\sl Proof.} Introduce elements $\Hb_{i,k}$ in $\Ub$,
$i\in[n]$ and $k\in\ZZ^\times$, such that 
$$\Kb_i^\pm(z)=\Kb_i^{\pm 1}\exp\left(\pm(q-q^{-1})\sum_{k\geq 1}
\Hb_{i,\pm k}\cdot z^{\mp k}\right).$$
Each $\Hb_{i,k}$ can be expressed as an algebraic expression in the
$\Kb_{i,k}$'s, more precisely
$$\sum_{k\geq 1}\Hb_{i,\pm k}z^{\mp k}=
\pm(q-q^{-1})^{-1}\log\biggr(\Kb_i^{\mp 1}\Kb_i^\pm(z)\biggl).$$
A direct computation shows that
in the component $\Rb^{(\vb)}$ of the polynomial representation
the element $\Hb_{i,\pm k}$, $k\in\NN^\times$, acts via multiplication 
by the polynomial $\hat\Hb_{i,\pm k}$ such that
$$\hat\Hb_{i,\pm k}=-{[k]\over k}\left(
t^{\mp k}\sum_{l=1}^{\overline v_{i-1}}x_l^{\pm k}
+t^{\pm k}\sum_{l=1+\overline v_i}^dx_l^{\pm k}\right).$$
We are done.
\Box

\vskip3mm

\noindent{\sl Proof of the theorem.}
The elements $\Ec_{i,k}$, $\Fc_{i,k}$ and $\Kc_{i,k}$
of $K^G(Z)_t$ decompose as the sum of their components in each
$K^G(Z_{\vb,\wb})_t$, where $\vb,\wb\in\Vb$.
Lemma 14 implies that all these components 
belong to $\Psi_t(\Ub_t)$. On the other hand Lemma 15 implies that
$$K^G(Z_{_A})_t\subset\Psi_t(\Ub_t)$$ 
for any diagonal matrix $A\in\Mb$.
The surjectivity now follows from Proposition 10. 
\Box

\vskip1cm

\noindent{\bf 3. The simple and standard modules.}

\vskip3mm

\noindent
We keep the previous notation. In particular $G=GL_d\times\CC^\times$,
$F=\coprod_\vb F_\vb$ is the (disconnected) $n$-step flag manifold,
$Z$ is the Steinberg-type variety, $M_\vb=T^*F_\vb$ the cotangent bundles, 
$M=\coprod_\vb M_\vb$, and
$$N=\{x\in{\rm End}(\CC^d)\,|\,x^n=0\}.$$
Moreover for any $x\in N$ let
$F_{\vb,x}=\{D\in F_\vb\,|\, x(D_i)\subset F_{i-1}\}$
be the Springer variety.
Let $T\i GL_d$ be the subgroup of diagonal matrices and set 
$A=T\times\CC^\times$. Then, the Bivariant Localization Theorem
[CG, Theorem 5.11.10] and the Bivariant Riemann-Roch Theorem 
[CG, Theorem 5.11.11] give a homomorphism of algebras
$$K^G(Z)\twoheadrightarrow H_*(Z^a)$$
for all $a=(s,t)\in A$, where $Z^a$ denotes the fixed point subvariety
and $H_*$ stands for the Borel-Moore homology with complex coefficients
(see either [GV]). Moreover this map is surjective [CG, Theorem 6.2.4].
As a consequence, for any $a=(s,t)\in A$
we have a surjective homomorphism
$$\Psi_a\,:\,\Ub_t\twoheadrightarrow K^G(Z)_t
\twoheadrightarrow H_*(Z^a).$$
Fix $a=(s,t)\in T\times\CC^\times$. 
Let $G(s)$ be the centralizer of $s$ in $GL_d$.
Let $M^a=\coprod_\vb M^a_\vb$ be the fixed point subvariety
and $\pi\,:\,M^a\to N^a$ the projection to
$$N^a=\{x\in N\,|\,s\,x\,s^{-1}=t^{-2}x\}.$$
The map $\pi$ commutes with the action of $G(s)$.
For all $\vb\in\Vb$ consider the $G(s)$-equivariant complex on $N^a$
$$\Lc_\vb^a=\pi_*\CC_{M^a_\vb}[\dim\,M^a_\vb]$$ 
and set $\Lc^a=\oplus_\vb\Lc^a_\vb$.
The equivariant version of the
Beilinson-Bernstein-Deligne Decomposition Theorem gives 
$$\Lc^a_\vb=\bigoplus_{O\subset N^a\atop i\in\ZZ}L_{O,\vb,i}\otimes IC_O[i],$$
where $IC_O[i]$ is the shift of the intersection cohomology complex
associated to the constant sheaf on $O$ (recall that
the isotropy subgroup in $G(s)$ of a point in $O$ is connected).
Put $L_{O,\vb}=\oplus_i L_{O,\vb,i}$ and $L_O=\oplus_\vb L_{O,\vb}$
for all $O$. The $L_O$'s are simple $\Ub_t$-modules
[GV, Theorem 8.4], [CG, Theorem 8.6.12] (for their non-vanishing see the remark 
after Theorem 4). The total cohomology space $H_*(F^s_x)$ 
is endowed with the structure of a $\Ub_t$-module
(see [CG, Proposition 8.6.15]).
Since this module only depends on the $G(s)$-orbit, $O$,
of $x$ in $N^a$, let denote it by $M_O$. The module $M_O$
is called a {\it standard module}. Then we have the following
analogue of the Kazhdan-Lusztig multiplicity formula 
(see [GV, Theorem 6.6] and [CG, Theorem 8.6.23]).

\vskip3mm

\proclaim Theorem 3. 
For any $a=(s,t)\in A$ and any $G(s)$-orbits $O,O'\in N^a$ one has 
$$[M_{O'}: L_O]=\sum_i\dim\Hc^i_{O'}(IC_O).$$
\Box

\vskip3mm

\noindent 
In the formula above $\dim\Hc^i_{O'}(IC_O)$ is the rank of the
$i$-th local cohomology sheaf of the complex $IC_O$
at any point of $O'$.
We recall that the quantized universal enveloping algebra
of $\widehat{\got{sl}}\,(n)$ (with trivial central charge), 
denoted by ${\bf U'}$, may be viewed as the 
$\KK$-subalgebra of ${\bf U}$ generated by Fourier coefficients
of the series (see [DF])
$$\Eb_i(z),\quad\Fb_i(z),\quad\Kb_{i+1}^\pm(t^iz)\Kb_i^\pm(t^iz)^{-1},
\qquad i=1,2,...,n-1.$$

\vskip3mm

\proclaim Lemma 16. The $\Ub_t$-modules $L_O$ are simple $\Ub'_t$-modules.

\vskip3mm

\noindent{\sl Proof.} Immediate since Fourier coefficients of
$\prod_{i=1}^n\Kb_i^\pm(zt^{2i})$
are central in $\Ub_t$.
\Box

\vskip3mm

\noindent Fix $t\in\CC^\times$, not a root of unity.
Then, simple $\Ub'_t$-modules are parametrized by their Drinfeld polynomials
(the proof is due to Chari and Pressley, see [CP2]).
Recall further that a nilpotent matrix $y\in{\rm End}(\CC^d)$ is labelled by 
the partition 
$\lambda(y)=(\lambda_1(y)\geq\lambda_2(y)\geq...)$ such that
$\lambda_i(y)$ is the length of the $i$-th Jordan block of $y$.
If $y^n=0$ the dual partition, 
$\lambda^\vee(y)=(\lambda^\vee_1(y)\geq\lambda^\vee_2(y)\geq\cdots 
\geq\lambda^\vee_n(y))$, is such that
$$\lambda^\vee_i(y)={\rm dim}{\rm Ker}(y^i)-{\rm dim}{\rm Ker}(y^{i-1}),
\qquad i=1,2,...,n.$$
In particular $\lambda^\vee(y)$ may be viewed as a dominant weight of
${\got{gl}}\,(n)$ in the obvious way.
Denote by $\lambda^\vee(O)$ the dual partition of a representative
of any nilpotent orbit $O$.

\vskip3mm

\proclaim Theorem 4. Fix $a=(s,t)\in A$ and fix a $G(s)$-orbit
$O\subset N^a$. Put $\lambda^\vee(O)=
(\lambda^\vee_1,\lambda^\vee_2,...,\lambda^\vee_n)$. Then,
the Drinfeld polynomials of the simple $\Ub'_t$-module $L_O$ are given by
$$P_i(z)=\prod_{\lambda^\vee_{i+1}<k\leq\lambda^\vee_i}
(z-t^{i-2}s_k^{-1}),\qquad i=1,2,...,n-1,$$
where the $s_k$'s are such that the matrix $s$ is conjugate to 
$$\bigoplus_{i=1}^{\lambda^\vee_1}s_iD(\lambda_i)
\quad{\sl with}\quad
D(k)=\sum_{i=0}^{k-1}t^{-2i}E_{ii}\quad\forall k\in\NN^\times.$$

\vskip1mm

\noindent{\bf Remark.} It follows that the $L_O$'s are non zero for all $O$.
Moreover this non-vanishing result is still valid if $t$ is a root of 
unity.

\vskip3mm

\noindent We first mention the following simple result of linear algebra
(see [CG, Remark 4.2.2]).

\vskip3mm

\proclaim Lemma 17. If $x\in N$ and $\vb\in\Vb$, then
$$F_{\vb,x}\neq\emptyset\quad\Rightarrow\quad\vb\leq \lambda^\vee(x),$$
where $\leq$ is the standard order on the weight lattice of
${\got{gl}}(n)$.
\Box

\vskip3mm

\noindent{\sl Proof of Theorem 4.} 
By construction we have
$$L_{O,\vb}=\{u\in L_O\,|\,\Kb_i\cdot u =t^{d-v_i}u,\quad
\forall i=1,2,...,n\}.$$
In other words, $L_{O,\vb}\i L_O$ is the subspace of weight
$\sum_{i=1}^{n-1}(v_i-v_{i+1})\omega_i$ where $\omega_i$
is the $i$-th fundamental weight of ${\got{sl}}(n)$.
Fix $x\in O$ and set $\lambda^\vee=\lambda^\vee(x)$.
Since $F_{\lambda^\vee,x}$ reduces to the single flag
$$0\subseteq{\rm Ker}(x)\subseteq {\rm Ker}(x^2)\subseteq
\cdots\subseteq{\rm Ker}(x^n)=\CC^d$$
the weight $\lambda^\vee$ has multiplicity one in the $\Ub'_t$-module.
Moreover, it follows from Lemma 17 that all
the weights of $M_O$ are less or equal to $\lambda^\vee$.
By construction
$$\matrix{
L_{O,\vb}\neq\{0\}&\Rightarrow O\subset\pi(M^a_\vb),\hfill\cr\cr
&\Rightarrow F_{\vb,x}^s\neq\emptyset.\hfill
}$$
Thus, the weights of $L_O$ are less or equal to $\lambda^\vee$ again.
On the other hand Theorem 3 implies that

\vskip2mm

\item{(i)} $[M_O: L_O]=1,$

\vskip2mm

\item{(ii)} $[M_O: L_{O'}]\neq 0\Rightarrow O\i\overline{O'}.$

\vskip2mm

\noindent Fix another orbit $O'\neq O$ containing $O$ in its closure.
Then $\lambda^\vee(O')<\lambda^\vee(O).$
So (i) and (ii) imply that
the module $L_O$ has highest weight $\lambda^\vee$ and that
$L_{O,\lambda^\vee}=H^*(F^s_{\lambda^\vee,x})$. Since 
$F^s_{\lambda^\vee,x}$ reduces to a single point, on $L_{O,\lambda^\vee}\i L_O$
we get (see 2.2)
$$\Kb_{i+1}^\pm(t^iz)\Kb_i^\pm(t^iz)^{-1}=
\prod_{\lambda^\vee_{i+1}<k\leq\lambda^\vee_i}
\theta_{-1}(t^{1-i}s_kz^{-1}),$$
and we are done.
\Box

\vskip3mm

\noindent Let $\omega_1,\omega_2,...,\omega_{n-1}$
be the fundamental weights of ${\got sl}_n$. For any $\alpha\in\CC^\times$ 
and any $i=1,2,...,n-1$, let $V(\omega_i)_\alpha$ denote
the fundamental representation of $\Ub'$, i.e.
the simple finite dimensional $\Ub'$-module whose
Drinfeld polynomials are
$$P_j(z)=(z-\alpha)^{\delta_{ji}},\quad\forall j=1,2,...,n-1.$$

\noindent 
For any partition 
$\lambda=(\lambda_1\geq\lambda_2\geq ...\geq \lambda_l\geq 0)$
of $d$
and any $\alpha=(\alpha_1,\alpha_2,...,\alpha_l)\in(\CC^\times)^l$, 
consider the block diagonal $d\times d$-matrices
$$x_\lambda=\bigoplus_{i=1}^lJ(\lambda_i)
\quad{\rm with}\quad 
J(k)=\sum_{i=1}^{k-1}E_{i,i+1}\quad\forall k\in\NN^\times,$$
and 
$$s_{\lambda,\alpha}=\bigoplus_{i=1}^l\alpha_iD(\lambda_i),$$
(see the previous theorem). In particular 
$s_{\lambda,\alpha}\,x_\lambda\, s_{\lambda,\alpha}^{-1}=t^2\,x_\lambda.$
Denote by $O_{\lambda,\alpha}$
the $G(s_{\lambda,\alpha})$-orbit of $x_\lambda$ in $N^{a_{\lambda,\alpha}}$,
where $a_{\lambda,\alpha}=(s_{\lambda,\alpha},t)$.

\vskip3mm

\noindent{\bf Proposition 18.} {\it If $\lambda_i<n$ for all $i$ 
then the classes of $M_{O_{\lambda,\alpha}}$ and 
$$V(\omega_{\lambda_1})_{\alpha_1^{-1}t^{\lambda_1-2}}\otimes
V(\omega_{\lambda_2})_{\alpha_2^{-1}t^{\lambda_2-2}}\otimes\cdots\otimes
V(\omega_{\lambda_l})_{\alpha_l^{-1}t^{\lambda_l-2}}$$
in the Grothendieck ring of finite dimensional
$\Ub'_t$-modules are equal.}

\vskip 3mm

\noindent{\it Proof.} First, recall the geometric
construction of the Schur functor given in [GRV, Theorem 6.8]
between $\Ub_t$-modules
and $\Hb_t$-modules, where $\Hb_t$ is the affine Hecke algebra 
of ${\got gl}_d$ with quantum parameter $t$. 
Suppose that $n\geq d$.
Let $F_{c}$ be the variety of complete flags in $\CC^d$ and let
$M_{c}=T^*F_{c}$ be the cotangent bundle. Given $a\in A$, let
$\pi\,:\,M^a_{c}\to N^a$ be the projection and consider the complex
$\Lc^a_{c}=\pi_*\CC_{M^a_{c}}[\dim M^a_{c}].$ Put
$$\Ub_a={\rm Ext}(\Lc^a,\Lc^a),\quad
\Hb_a={\rm Ext}(\Lc_{c}^a,\Lc_{c}^a)
\qquad{\rm and}\qquad
\Wb_a={\rm Ext}(\Lc_{c}^a,\Lc^a)$$
(the Yoneda Ext-groups in the bounded derived category of constructible 
sheaves on $N^a$.)
The space $\Wb_a$ is a $(\Ub_a,\Hb_a)$-bimodule. Moreover,
since we have surjective algebra homomorphisms
$$\Ub_t\twoheadrightarrow\Ub_a
\qquad{\rm and}\qquad
\Hb_t\twoheadrightarrow\Hb_a,$$
a $\Hb_a$-module (resp. a $\Ub_a$-module) can be pulled back 
to a $\Hb_t$-module (resp. a $\Ub_t$-module).
The Schur dual of the pull-back of a $\Hb_a$-module $V$ is the pull-back
of the $\Ub_a$-module $V^o=\Wb_a\otimes_{\Hb_a}V.$
In particular, fix a nilpotent orbit $O\subset N^a$ and fix a point $x\in O$.
Denote by $j$ the embedding $\{x\}\hookrightarrow N^a$.
The standard module $M_O$ is isomorphic to the cohomology space
$H^*(j^!\Lc^a)$, where $\Ub_t$ acts via the algebra homomorphism
$$\Ub_t\twoheadrightarrow\Ub_a\to {\rm Ext}(j^!\Lc^a,j^!\Lc^a).$$
$M_O$ splits as a vector space into the direct sum
$$M_O=\bigoplus_{O',k}L_{O'}\otimes H^k(j^!IC_{O'}).$$
Similarly the underlying vector space of the
standard $\Hb_t$-module corresponding to $O$ is
$$N_O=\bigoplus_{O',k}L_{c,O'}\otimes H^k(j^!IC_{O'}).$$
Let $\Ab_a\subset\Hb_a$ be the maximal semi-simple subalgebra of $\Hb_a$, i.e.
$\Ab_a=\bigoplus_O{\rm End}(L_{c,O}).$
Recall that we have an isomorphism 
$\Wb_a=\biggl(\bigoplus_O{\rm Hom}(L_{c,O},L_O)\biggr)\otimes_{\Ab_a}\Hb_a$
(see [GRV; (6.12)]). Then,
$$\matrix{
N_O^o
&=\bigoplus_{O'',k}\biggl(\bigoplus_{O'}{\rm Hom}(L_{c,O'},L_{O'})\biggr)
\otimes_{\Ab_a}\Hb_a\otimes_{\Hb_a}L_{c,O''}\otimes H^k(j^!IC_{O''}),
\hfill\cr\cr
&=\bigoplus_{O',k}\biggl({\rm Hom}(L_{c,O'},L_{O'})\otimes_{\Ab_a}
L_{c,O'}\otimes H^k(j^!IC_{O'})\biggr),\hfill\cr\cr
&=M_O.\hfill
}$$
Now using the Induction Theorem [KL, Theorem 6.2] one proves that
the class $[N_O]$ of $N_O$ in the Grothendieck ring
is the same as the class of an induced module (see [A]).
More precisely, we have
$$[N_{O_{\lambda,\alpha}}]=[\Hb_t\otimes_{\Hb_{\lambda,t}}\CC_\alpha].$$
Here $\Hb_{\lambda,t}\subset\Hb_t$ is the affine
Hecke algebra of the Young subgroup of ${\got S}_d$
associated to the partition $\lambda$ and $\CC_\alpha$
is the one-dimensional module such that
$$T_i\mapsto t
\quad{\rm and}\quad
X_i\mapsto s_{\lambda,\alpha,i},$$
where $s_{\lambda,\alpha,i}$ is the $i$-th entry in $s_{\lambda,\alpha}$.
In particular if $\lambda=(d)$ then
$x_\lambda$ is regular and, thus,
$N_{O_{\lambda,\alpha}}$ is the simple $\Hb_t$-module labelled 
by $O_{\lambda,\alpha}$. If moreover $d<n$, its Schur dual is the
simple $\Ub_t$-module
$L_{O_{\lambda,\alpha}}\simeq V(\omega_d)_{\alpha^{-1}t^{d-2}}$.
Using the behaviour of the Schur functor with respect to tensor products 
(see [CP1] for the case of the affine Hecke algebra) we get
the result. 
\Box

\vskip 3mm

\noindent{\bf Remark.} 
Theorem 3 implies in particular that if 
$O_{\lambda,\alpha}\subset N^{a_{\lambda,\alpha}}$ is open then
$V(\omega_{\lambda_1})_{\alpha_1^{-1}t^{\lambda_1-2}}\otimes
V(\omega_{\lambda_2})_{\alpha_2^{-1}t^{\lambda_2-2}}\otimes\cdots\otimes
V(\omega_{\lambda_l})_{\alpha_l^{-1}t^{\lambda_l-2}}$
is irreducible.

\vskip1cm

\centerline{\bf References}

\vskip3mm

\hangindent3cm[A]\qquad\qquad Ariki, S.: 
On the decomposition number of the Hecke algebra of $G(m,1,n)$.
{\sl preprint}, (1996).

\vskip1mm

\hangindent3cm[BLM]\quad\quad Beilinson, A., Lusztig, G., MacPherson, R.:
A geometric setting for quantum groups.
{\sl Duke Math. J.}, {\bf 61} (1990), 655-675.

\vskip1mm

\hangindent3cm[CG]\quad\qquad Chriss, N., Ginzburg, V.:
Representation theory and complex geometry.
{\sl Birkhauser} (1997).

\vskip1mm

\hangindent3cm[CP1]\quad\qquad Chari, V., Pressley, A.:
Quantum affine algebras and affine Hecke algebras.
{\sl qalg-preprint}, {\bf 9501003}.

\vskip1mm

\hangindent3cm[CP2]\quad\qquad Chari, V., Pressley, A.:
A guide to quantum groups.
{\sl Cambridge University Press} (1994).

\vskip1mm

\hangindent3cm[D]\qquad\qquad Drinfeld, V.: 
A new realization of Yangians and quantized affine algebras.
{\sl Soviet. Math. Dokl.}, {\bf 36} (1988), 212-216.

\vskip1mm

\hangindent3cm[DF]\quad\qquad Ding, J., Frenkel, I.:
Isomorphism of two realizations of the quantum affine
algebra $U_q(\widehat{\got {gl}}(n))$.
{\sl Comm. Math. Phys.}, {\bf 156} (1994), 277-300.

\vskip1mm

\hangindent3cm[GRV]\quad\quad Ginzburg, V.,  Reshetikhin, N.,  Vasserot, E.:
Quantum groups and flag varieties.  {\sl Contemp. Math.}, {\bf 175} (1994),
101-130.

\vskip1mm

\hangindent3cm[GV]\quad\qquad Ginzburg, V., Vasserot, E.:
Langlands reciprocity for affine quantum groups of type $A_n$.
{\sl Internat. Math. Res. Notices}, {\bf 3} (1993), 67-85.

\vskip1mm

\hangindent3cm[KL]\quad\qquad Kazhdan, D., Lusztig, G.:
Proof of the Deligne-Langlands conjecture for Hecke algebras.
{\sl Invent. Math.}, {\bf 87} (1987), 153-215.

\vskip1mm

\hangindent3cm[M]\quad\qquad\  MacDonald, I.G.: 
Symmetric functions and Hall polynomials.
{\sl Oxford university press} (1995).

\vskip1mm

\hangindent3cm[V]\qquad\qquad  Vasserot, E.: 
Representations de groupes quantiques et permutations.
{\sl Annales Sci. ENS}, {\bf 26} (1993), 747-773.

\vskip3cm

$$\matrix{
{\rm Eric\ Vasserot}\hfill\cr
{\rm D\acute{e}partement\ de\ Math\acute{e}matiques}\hfill\cr
{\rm Universit\acute{e}\ de\ Cergy-Pontoise}\hfill\cr
{\rm 2\ Av.\ A.\ Chauvin}\hfill\cr
{\rm 95302\ Cergy-Pontoise\ Cedex}\hfill\cr
{\rm France}\hfill\cr
{\rm email:\ vasserot\ at\ math.pst.u-cergy.fr}\hfill
}$$

\bye